\documentclass{article}

\usepackage{cmd}
\usepackage{format}

\title{Marginalia to a Theorem of Asper\'o and Schindler}
\author{Obrad Kasum\footnote{Mr.\ Kasum has received funding from the European Union’s Horizon 2020 research and innovation program under the Marie Skłodowska-Curie grant agreement No.\ 945322 \includegraphics[height=2.5mm]{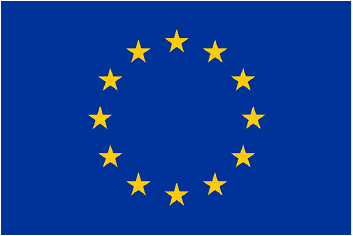}}, Boban Veli\v ckovi\' c}
\date{December 18, 2024}

\begin{document}

\maketitle
\begin{abstract}
    We give a game-theoretic characterization of when a model of an infinitary propositional formula can be added by a proper, semiproper, and stationary-set-preserving poset.
    In the latter case, we also give a general sufficient condition for the existence of such a poset.
    We use this condition to give a somewhat different proof of the theorem of Asper\'o and Schindler, which states that $\mathsf{MM}^{++}$ implies Woodin's axiom $(*)$.
\end{abstract}
\tableofcontents

\section{Introduction}

One of the fundamental problems in logic is to characterize when a given formula has a model and what that model is.
The most important instance of this problem is the well-known \textit{Entscheidungsproblem}, which asks whether a first-order theory is satisfiable.
In this paper, we will be interested in certain forms of the problem which concern infinitary propositional formulas.
In this case, the most basic form of a model-existence property is, of course, the \textit{consistency}.

\begin{definition}
    Let $\phi$ be an infinitary propositional formula.
    Then $\phi$ is \intro{consistent} iff there exists a poset $\P$ such that in $V^\P$, there exists a model\footnote{A model for $\phi$ is an evaluation for its propositional letters which assigns the value $1$ to $\phi$.} for $\phi$.
    \qed
\end{definition}

A superficial inspection of this definition might cause doubts as to whether it is natural to look at models in all generic extensions.
These doubts are quickly eliminated once we realize that this notion is equivalent to the property of ``not proving the contradiction''.
To fully appriciate this last comment, one has to recall the fact that the infinitary propositional logic comes equipped with a natural notion of a proof system, one which is a relatively straightforward generalization of the proof system for the standard (finitary) propositional logic.\footnote{See, for example, \cite[§1]{farah2021extender}.}
This explains the central role in Logic of this notion of consistency.
The starting point of our investigation is the question of what happens with this property if we restrict the class of posets which are to be considered for adding a model to the formula.

\begin{definition}
    Let $\phi$ be an infinitary propositional formula.
    Then $\phi$ is \intro{proper-consistent} (resp. \intro{semiproper-consistent}, \intro{ssp-consistent}\footnote{The abbreviation ``\intro{ssp}'' stands for ``stationary set preserving''.}) iff there exists a proper (resp. semiproper, ssp) poset $\P$ such that in $V^\P$, there exists a model for $\phi$.
    \qed
\end{definition}

For each of these three notions, we will give a game-theoretic characterization.
The motivation for such pursuit is the fact that the notion of consistency has, in addition to the proof-theoretic characterization, an important game theoretic one.
The latter is given in terms of the \textit{Model Existence Game}\footnote{See \cite[Section 8.3]{vaananen2011models}.
The game presented there greatly simplifies when we are only interested in propositional formulas.}.

\begin{definition}\label{meg}
    Suppose that $\phi$ is an infinitary propositional formula in negation normal form\footnote{A formula is a \intro{negation normal form} iff all the negations are immediately in front of propositional letters. Every formula is canonically equivalent to a formula in negation normal form. This form is obtained simply by pushing the negations ``all the way inside'' using De Morgan's laws.}.
    Then the game \intro{$\Game^\mathrm{meg}(\phi)$} is defined as the length $\omega$ two player game of the form
    \begin{center}
        \begin{tabular}{c|cccccc}
            I  &      $Q_0$ &       & $Q_1$ &       & $\cdots$\\
            \hline
            II &       & $w_0$ &       & $w_1$ & $\cdots$
        \end{tabular}
    \end{center}
    where $w_{-1}:=\{\phi\}$ and for all $n<\omega$, the following is satisfied.
    \begin{parts}
        \item Player I must ensure that either
        \begin{parts}
            \item $Q_n=\psi$ for some $\psi\in w_{n-1}$ of the form $\bigvee_{i\in I}\theta_i$, or

            \item $Q_n=(\psi,i)$ for some $\psi\in w_{n-1}$ of the form $\bigwedge_{i\in I}\theta_i$ and some $i\in I$.
        \end{parts}

        \item Player II must ensure that if $Q_n=\psi$ for some $\psi\equiv\bigvee_{i\in I}\theta_i$, then for some $i\in I$, $\theta_i\in w_n$.

        \item Player II must ensure that if $Q_n=(\psi,i)$ for some $\psi\equiv\bigwedge_{i\in I}\theta_i$, then $\theta_i\in w_n$.
    \end{parts}
    The infinite plays with no rules broken are won by Player II.
    \qed
\end{definition}

The characterization of consistency that we want to imitate reads as follows.

\begin{theorem}
    Let $\phi$ be an infinitary propositional formula in negation normal form.
    Then $\phi$ is consistent if and only if Player II has a winning strategy in $\Game^\mathrm{meg}(\phi)$.
    \footnote{
    This follows by the proof of \cite[Theorem 8.12]{vaananen2011models}.
    The point is that once the set of all subformulas of $\phi$ is made countable, we are in the scenario of the quoted theorem.
    Since the game is closed for Player II, it is generically absolute whether or not Player II wins.
    }
    \qed
\end{theorem}

We will describe versions of the Model Existence Game that correspond to the notions of proper-, semi-proper, and ssp-consistency.
They are denoted by $\Game^\mathrm{p}_\kappa(o)$, $\Game^\mathrm{sp}_\kappa(o)$, and $\Game^\mathrm{ssp}_\kappa(o)$ (respectively), and depend implicitly on $\phi$ and explicitly on a beth-fixed point $\kappa$.
We prove the following.\footnote{See Corollaries \ref{878}, \ref{1702}, and \ref{1993}.}

\begin{theorem}
    Let $\phi$ be an infinitary propositional formula in negation normal form.
    Suppose there exists a proper class of inaccessible cardinals.
    Then the following are equivalent.
    \begin{parts}
        \item $\phi$ is proper- (resp. semiproper-, ssp-) consistent.

        \item For some inaccessible $\kappa$, Player II has a winning strategy in $\Game^\mathrm{p}_\kappa(o)$ (resp. $\Game^\mathrm{sp}_\kappa(o)$, $\Game^\mathrm{ssp}_\kappa(o)$).
        \qed
    \end{parts}
\end{theorem}

This theorem essentially characterizes when a model to a $\boldsymbol{\Sigma}_1$ formula can be added by forcing of the given kind.
Let us mention some examples.
There is a formula whose models (in all generic extensions) naturally code surjections $\omega \twoheadrightarrow \omega_1^V$.
This is an example of a consistent formula which is not ssp-consistent.
On the other hand, we will see a formula whose models naturally code a cofinal map $\omega \to \omega_2^V$.\footnote{See Example \ref{Namba problem}.}
This is an example of an ssp-consistent formula which is not proper-consistent.\footnote{The fact that the formula is ssp-consistent follows from the fact that the Namba forcing is always ssp.
We prove it explicitly in Example \ref{2398}, but by assuming more than $\ZFC$.
On the other hand, it is a matter of routine to verify that the cofinality of $\omega_2^V$ cannot be changed to $\omega$ by a proper poset, which explains why the formula is not proper consistent.
}
An example of a proper-consistent formula drawn from the literature is given in Section \ref{An Example for the Proper Consistency}.\footnote{The formula is given implicitly, since we work with the notion of a \textit{problem}.
This notion and its exact relationship to infinitary propositional formulas is described in Section \ref{The Notion of a Problem}.}

In addition to simply answering the question of whether a model for a given formula can be forced by a poset of a certain kind, our method provides such a poset whenever the answer is positive.
The poset thus obtained is made up of conditions consisting of working parts and \textit{side-conditions}.
This is exactly the form of proper posets that dominated the literature on applications of $\PFA$.\footnote{See, for example, \cite{abraham2007some} or \cite{moore2006five}.}
Thus, our result can also be seen as an abstract unification theorem for those seemingly ad hoc constructions.

Finally, in a recent breakthrough of D. Asperó and R. Schindler,\footnote{\cite{aspero2021martin}} a new sophisticated construction of ssp posets was presented.
We observe that the concrete aim towards which their poset was built can be expressed by an infinitary propositional formula.
One is then confronted by the question of whether that formula is ssp-consistent.
Once this is verified, our method provides an ssp poset forcing a model for that formula.\footnote{
Let us note here that the Asperó-Schindler poset is somewhat analogous to the posets made up of pairs of working parts and side conditions, but it did not really fit into this paradigm, while our poset is exactly of this form.}
Thus, the final question that we want to answer is how to verify the ssp-consistency for the formulas of ``Asperó-Schindler type''.\footnote{For example, in addition to the construction of \cite{aspero2021martin} itself, its variation plays an important role in the argument of \cite{lietz2023forcing}.}

We resolve this final question by providing an abstract condition on infinitary propositional formulas that we call the \textit{AS-goodness} and show that it implies ssp-consistency.
This condition is stated in terms of generic elementary embeddings and we will verify that the concrete formula relevant for the Asperó-Schindler result satisfies it.
This verification is also an explanation of why we expect this condition to hold in other applications of their method\footnote{Such as, for example, the one that appears in \cite{lietz2023forcing}.}.

Let us give a brief outline of the paper.
\begin{description}
    \item[Sec. \ref{The Notion of a Problem}] We define a simply combinatorial notion of a \textit{problem} and explain why it is same whether we work with infinitary propositional formulas or with problems.
    We will then, for the rest of the paper, mostly work with problems, which is essentially just a matter of simplifying some proofs.

    \item[Sec. \ref{Some Simple Remarks}] We recall the characterization of genericity and semigenericity which will play the main role in designing our game-theoretic characterizations.
    Additionally, we introduce the notion of a virtual model and some basic notation around it.

    \item[Sec. \ref{Proper Consistency}] This is where we present the game-theoretic characterization of proper consistency.

    \item[Sec. \ref{An Example for the Proper Consistency}] Here we give an example a formula whose proper consistency has been (implicitly) considered in the literature.
    We illustrate how to apply our method in this case.

    \item[Sec. \ref{Semiproper Consistency}] This is where we present the variation of our method that concerns the semi-proper consistency.

    \item[Sec. \ref{ssp consistency}] Here we give the variation of the method for the ssp consistency.

    \item[Sec. \ref{AS Goodness}] In this section, we prove the ssp-consistency criterion based on generic embeddings.

    \item[Sec. \ref{An Example for the AS Goodness}] Finally, we illustrate that criterion on the example of a formula coding the main task of \cite{aspero2021martin}.
\end{description}


\section{The Notion of a Problem}
\label{The Notion of a Problem}

Given an infinitary propositional formula, we can think of it as a problem.
A solution to this problem is a model for the formula, i.e. an evaluation that assigns the value $1$ to the formula.
The following poset is a natural attempt to force a model of a formula.

\begin{definition}
    Let $\phi$ be a consistent infinitary propositional formula in negation normal form.
    Then the poset \intro{$\Hi(\phi)$} consists of all finite consistent sets of subformulas of $\phi$ that contain $\phi$, ordered by the reverse inclusion.
    \qed
\end{definition}

Forcing with this poset adds a model of the formula, but it will usually collapse $\omega_1$.
If we want to force a solution with, say, a proper forcing, then we must not produce generics for $\Hi(\phi)$.
Luckily, we do not actually need full generics.

\begin{definition}
    Let $\phi$ be a consistent infinitary propositional formula in negation normal form.
    Then the set \intro{$\dd(\phi)$} of dense subsets of $\Hi(\phi)$ is defined to be
    $$\{D_{\psi},\, D_{\theta,i}\}_{\psi,\theta,i},$$
    where
    \begin{parts}
        \item for $\psi$ a subformula of $\phi$ with $\psi\equiv\bigvee_{i\in I}\phi_i$,
        $$D_\psi:=\{w\in\Hi(\phi) : \psi\in w\implies\exists i\in I, \psi_i\in w\},$$

        \item for $\theta$ a subformula of $\phi$ with $\theta\equiv\bigwedge_{i\in I}\theta_i$ and for $i\in I$,
        $$D_{\theta,i}:=\{w\in\Hi(\phi) : \theta\in w\implies \theta_i\in w\}.$$
        \qed
    \end{parts}
\end{definition}

It turns out that adding a filter for $\Hi(\phi)$ which is just $\dd(\phi)$-generic\footnote{In other words, a filter that intersects all the sets in $\dd(\phi)$, but not necessarily all dense sets in $V$.} suffices.

\begin{proposition}
    Let $\phi$ be a consistent infinitary propositional formula in negation normal form.
    Then there exists a $\dd(\phi)$-generic filter for $\Hi(\phi)$ if and only if there exists a model for $\phi$.
\end{proposition}
\begin{proof}
    If we are given a $\dd(\phi)$-generic filter $g$ for $\Hi(\phi)$, we can set
    $$\mu_g(l):=\begin{cases}
        1, & \exists w\in g, l\in w\\
        0, & \mbox{otherwise,}
    \end{cases}$$
    for all propositional letters $l$ appearing in $\phi$.
    It then follows by the induction on the complexity of a subformula that $\mu_g\models\phi$.
    On the other hand, if we are given a model $\mu\models\phi$, we define
    $$g_\mu:=\{w\in\Hi(\phi) : \mu\models w\}$$
    and routinely verify that $g_\mu$ is a $\dd(\phi)$-generic filter for $\Hi(\phi)$.
\end{proof}

Thus, we are lead to the following generalization of the notion of a problem given by a propositional formula.

\begin{definition}
    A pair $(\Hi,\dd)$ is a \intro{problem} iff
    \begin{parts}
        \item $\Hi$ is a poset\footnote{For the sake of precision, we will mean by ``poset'' a partial order with the distinguished maximum.},

        \item for all $p,q\in\Hi$, if $p\para q$, then $\inf\{p,q\}$ exists (and it will be denoted by \intro{$pq$}),

        \item $\dd$ is a collection of pre-dense subsets of $\Hi$.
        \qed
    \end{parts}
\end{definition}

A solution for a problem then becomes simply a $\dd$-generic filter for $\Hi$.
The proper consistency then naturally becomes the property of being able to force a solution by a proper forcing.

\begin{definition}
    A problem $(\Hi,\dd)$ is \intro{proper-consistent} (resp. \intro{semiproper-consistent}, \intro{ssp-consistent}) iff there exists a proper (resp. semiproper, ssp) poset $\Q$ such that in $V^\Q$, there exists a $\dd$-generic filter for $\Hi$.
    \qed
\end{definition}

\begin{lemma}
    Let $\phi$ be a consistent infinitary propositional formula in negation normal form.
    Then $\phi$ is proper- (resp. semiproper-, ssp-) consistent if and only if $(\Hi(\phi),\dd(\phi))$ is so.
    \qed
\end{lemma}

Not all problems come from an infinitary propositional formula, but the most complex ones do.
The possibilities of the expressing power of propositional formulas will be illustrated in Section \ref{An Example for the AS Goodness}.
The apparent increased generality that problems introduce is just that, apparent.
We explain further in the following example.

\begin{example}
    Let $(\Hi,\dd)$ be a problem.
    We want to find a propositional formula $\chi$ such that there is a canonical correspondence between $\dd$-generic filters for $\Hi$ and models for $\chi$ (in all generic extensions).
    We will use the set
    $$\{g_w : w\in\Hi\}$$
    of propositional letters to build $\chi$.
    The idea is that the evaluation for $g_w$ answers ``yes'' or ``no'' to the question of whether $w$ belongs to the corresponding filter.
    The formula $\chi$ is then the conjunction of the following formulas:
    \begin{parts}
        \item $g_{1_\Hi}$,

        \item $\bigwedge_{w\leq w'}(g_w\implies g_{w'})$,

        \item $\bigwedge_{w,w'\in\Hi}(g_w\wedge g_{w'}\implies g_{ww'})$,

        \item $\bigwedge_{D\in\dd}\bigvee_{w\in D}g_w$.
    \end{parts}
    
    Now, in any generic extension, if $g$ is a $\dd$-generic filter for $\Hi$, then
    $$\{g_w : w\in\Hi\}\to 2:g_w\mapsto\begin{cases}
        1, & w\in g\\
        0, & w\not\in g
    \end{cases}$$
    is a model for $\chi$.
    On the other hand, if $\mu$ is a model for $\chi$, then
    $$\{w\in\Hi : \mu(g_w)=1\}$$
    is a $\dd$-generic filter for $\Hi$.
    Furthermore, these two correspondences are mutual inverses.
    \qed
\end{example}

The next example illustrates that sometimes it might be easier to work with the abstract notion of a problem than to pass through a propositional formula.

\begin{example}\label{Namba problem}
    We explain how to code the ``Namba problem''.
    More precisely, we want to force a cofinal function $\omega\to\omega_2^V$ by an ssp poset.
    Such a function is a subset of $\omega\times\omega_2$, so we can build a natural propositional formula $\phi$ with the propositional letters
    $$\{f_{n,\alpha} : n<\omega,\alpha<\omega_2\}$$
    such that for every evaluation $\nu$ in every generic extension, $\nu\models\phi$ if and only if the set
    $$\{(n,\alpha) : \nu(f_{n,\alpha})=1\}$$
    is a cofinal function $\omega\to\omega_2^V$.
    The Namba problem then becomes the question whether $\phi$ is ssp-consistent, or equivalently, whether $(\Hi(\phi),\dd(\phi))$ is ssp-consistent.

    On the other hand, the generality that we get from the notion of a problem allows us to code the Namba problem more directly, i.e. without passing through a propositional formula.
    We then set
    $$\Hi:=\{w\subseteq\omega\times\omega_2 : w\mbox{ is a finite function}\},$$
    ordered by the reverse inclusion, and we set $\dd:=\{D_n,C_\alpha : n<\omega,\alpha,\omega_2\}$, where
    $$D_n:=\{w\in\Hi : n\in\dom(w)\},$$
    $$C_\alpha:=\{w\in\Hi : \exists m\in\dom(w),w(m)>\alpha\}.$$
    Then a filter $F$ for $\Hi$ (in a generic extension) is $\dd$-generic if and only if 
    $$\bigcup F:\omega\xrightarrow[]{\mathrm{cof}}\omega_2^V.$$
    \qed
\end{example}


\section{Some Simple Remarks}
\label{Some Simple Remarks}

In our arguments, we will define a class of preconditions, and then we will define a way to pick among them those that should be conditions in our final posets.
The preconditions will simply be pairs of a working part together with a chain of side conditions.
The criterion for picking the conditions is based on games.
The games are designed to ensure that the final poset adds a correct object.
However, they must also ensure that a condition of the poset is (semi)generic for a model that appears as a side condition.
To do this, we use the following slight modification of the standard (semi)genericity criterion.

\begin{lemma}\label{94}
    Suppose that
    \begin{assume}
        \item $\Q$ is a poset,

        \item $\theta\gg\rank(\Q)$ is regular,

        \item $M\prec (H_\theta,\in,\Q)$ is countable,

        \item $p\in\Q$.
    \end{assume}
    Then the following are equivalent.
    \begin{parts}
        \item\label{187} $p$ is generic (resp. semigeneric) for $(M,\Q)$.

        \item\label{189} For all $q\leq p$, for all $E\in M$ which contain $q$, there exists $r\in E\cap\Q$ such that $r\in M$ (resp. $\delta(\hull(M,r))=\delta(M)$)\footnote{We write \intro{$\delta(M)$} or \intro{$\delta_M$} for the intersection of $M$ with $\omega_1$. We use the same notation for all \textit{virtual models} (cf. Definition \ref{311}).} and $r\parallel q$.
    \end{parts}
\end{lemma}
\begin{proof}
    The modification here is that $E$ is not required to be (pre)dense.
    In order to apply the standard lemma, we just need to replace $E$ by $E\cup E^\bot$, where
    $$E^\bot:=\{r\in\Q : \forall q'\in E,r\bot q'\}.$$
\end{proof}

The side conditions in our posets will be \textit{virtual models}.
They were first defined by Veli\v ckovi\'c when studying alternative ways to iterate semiproper posets, see for example \cite{kasum2023iterating}.
In the cited paper, one can also find a detailed elaboration on how to work with virtual models.
We review here only the basics.
A prototype of a virtual model is an elementary submodel of some $H_\theta$, with $\theta$ regular, but there is more to the notion.

\begin{definition}\label{311}
    A set $M$ is a \intro{virtual model} iff $M\prec\widehat{M}\models\ZFC^-$.\footnote{\intro{$\widehat M$} denotes the transitive closure of $M$.}
    \qed
\end{definition}

One of the fundamental operations for virtual models is the hull operation.

\begin{definition}
    Let $M$ be a virtual model and $X$ an arbitrary set.
    Then \intro{$\Hull(M,X)$} is defined as the set of all $f(x)$ where $f$ is a function in $M$ and $x$ is an element of $X\cap\dom(f)$.
    \qed
\end{definition}

$\Hull(M,X)$ is elementary in $\widehat M$ and it contains $X\cap\widehat M$.
In fact, it is the smallest virtual model with respect to these two properties.
The most important virtual models which are \underline{not} elementary submodels of some $H_\theta$ originate from an application of the following operation.

\begin{definition}
    Suppose that $M$ is a virtual model and that $\lambda$ is an infinite ordinal.
    We define \intro{$M\proj\lambda$} to be the image of $M$ under the transitive collapse of $\Hull(M,V_\lambda)$.
    \qed
\end{definition}

The anti-collapse of a virtual model is essentially a partial extender.
The following parameter measures the strength of that extender.

\begin{definition}
    Suppose that $M$ is a virtual model.
    Then \intro{$\lambda_M$} is defined as the largest $\lambda$ such that $V_\lambda\in M$, when it exists.
    \qed
\end{definition}

We will mostly talk about virtual models which belong to $\Cor_\lambda$, for a beth-fixed point $\lambda$.
These collections are defined just below.
Since we will refer a lot to beth-fixed points, it helps to introduce the following notation.

\begin{notation}
    \intro{$\Bf$} denotes the class of all beth-fixed points.
    \qed
\end{notation}

\begin{definition}
    Suppose that $\lambda\in\Bf$.
    Then \intro{$\Cor_\lambda$} denotes the set of all virtual model $M$ satisfying that
    \begin{parts}
        \item $\lambda_M=\lambda$,

        \item $\widehat{M}=\hull(M,V_{\lambda})$.
        \qed
    \end{parts}
\end{definition}

\begin{definition}
    The class \intro{$\Cor$} is defined as $\bigcup_{\lambda\in\Bf}\Cor_\lambda$.
    \qed
\end{definition}

If $M\in\Cor$ and if $\Q$ is a poset that belongs to $M\cap V_{\lambda_M}$, we can make sense of the phrase ``$\Q$ is (semi)proper for $M$''.
One way to do so is to just run the usual definition, but with respect to this more general object $M$.\footnote{
The point is that $M$ is not necessarily an elementary submodel of some $H_\theta$.
}
It turns out that this notion is equivalent to saying that $\Q$ is (semi)proper for $M\cap H_\theta$ for any and all regular $\theta\in (2^{2^\kappa},\lambda_M)\cap M$.
However, $M\cap H_\theta\prec H_\theta$, so we do not diverge much from the usual notion of (semi)properness.

When using virtual models as side conditions, we will require them to be $\in$-chains.
In fact, we will require a bit more, as explained in the following definition.

\begin{definition}
    Suppose that $\M\subseteq\Cor$.
    Then $\M$ is a \intro{vm-chain} iff for all $M,N\in\M$,
    \begin{parts}
        \item if $\delta_M=\delta_N$, then $M=N$,
        
        \item if $\delta_M<\delta_N$, then $M\in N$ and $\lambda_M<\lambda_N$.
        \qed
    \end{parts}
\end{definition}

When passing from a virtual model $M$ to its projection $M\proj\lambda$, we lose information.
In a sense, $M$ then certifies that $M\proj\lambda$ is ``consistent'' with more information.\footnote{
Our intuition is that it is more complicated to embed a model into a transitive structure of a larger rank.
}
This is not true for all models $N\in\Cor_\lambda$ and even if it is true for some such $N$, a virtual model $M$ such that $M\proj\lambda=N$ is not necessarily unique.
We will use the following notion in order to talk about this possibility of increasing the information content of a virtual model.

\begin{definition}
    Suppose that
    \begin{assume}
        \item $M$ and $N$ are virtual models,
        
        \item $M\in\Cor$,

        \item $\pi:M\longrightarrow N$.
    \end{assume}
    Then $\pi$ is a \intro{lifting} iff, letting
    $$\rho:\widehat{N\proj\lambda_M}\xrightarrow{\ \simeq\ }\hull(N,V_{\lambda_M})\prec\widehat{N}$$
    be the anti-collapse, we have that $M=N\proj\lambda_M$ and $\pi=\rho\rest M$.
    \qed
\end{definition}


\section{Proper Consistency}
\label{Proper Consistency}

\begin{declaration}
    We fix a problem $(\Hi,\dd)$.
    \qed
\end{declaration}

We want to characterize when this problem is proper-consistent.
In this case, we will also give a poset that adds a $\dd$-generic filter for $\Hi$.
We start by defining a class-sized poset of preconditions.
The required will be a (set-sized) suborder of it. 

\begin{definition}
    \hfill
    \begin{parts}
        \item The class \intro{$\P^\mathrm{p}$} consists of all $p=(w_p,\M_p)$ such that $w_p\in\Hi$ and $\M_p$ is a finite vm-chain.

        \item The order \intro{$\leq^\mathrm{p}$} on vm-chains is the reverse inclusion $\supseteq$.

        \item The order \intro{$\leq^\mathrm{p}$} on $\P^\mathrm{p}$ is defined by asserting that $p\leq q$ iff $w_p\leq_\Hi w_q$ and $\M_p\leq^\mathrm{p}\M_q$.
        \qed
    \end{parts}
\end{definition}

The criterion for deciding which preconditions should be taken as conditions in the final poset is based on a game.
We describe this game in the following definition.
Player I plays by asking certain questions and Player II plays by answering those questions.
Preconditions will be conditions if Player II (the ``good'' player) can always play the corresponding game for $\omega$ rounds, i.e. without failing to produce an answer of the required form.
There are three types of questions that Player I can ask.
The first ensures that the final poset actually forces a $\dd$-generic filter for $\Hi$.
The other two types of questions ensure that the final poset is proper.
They are based on the genericity criterion given in Lemma \ref{94}.

\begin{definition}
    Suppose that $\kappa\in\Bf$ and $p\in\P^\mathrm{p}\cap V_\kappa$.
    Then the game \intro{$\Game^\mathrm{p}_\kappa(p)$} is defined as the length $\omega$ two Player game of the form
    \begin{center}
        \begin{tabular}{c|cccccc}
            I  &       & $Q_0$ &       & $Q_1$ &       & $\cdots$\\
            \hline
            II & $p_{-1}$ &       & $p_0$ &       & $p_1$ & $\cdots$
        \end{tabular}
    \end{center}
    where $p_{-1}\in\P^\mathrm{p}\cap V_\kappa$ is such that $p_{-1}\leq p$ and for all $n<\omega$, the following is satisfied.
    \begin{parts}
        \item Player I must ensure that either
        \begin{parts}
            \item $Q_n=D$ for some $D\in\dd$, or

            \item $Q_n=M$ for some countable $M\prec H((2^\kappa)^+)$, or

            \item $Q_n=(M,E)$ for some $M\in\M_{p_{n-1}}$ and some $E\in M$.
        \end{parts}

        \item Player II must ensure that $p_n\in\P^\mathrm{p}\cap V_\kappa$ and $p_n\leq p_{n-1}$.

        \item Player II must ensure that if $Q_n=D$ for some $D\in\dd$, then there exists $w\in D$ such that $w_{p_n}\leq w$.

        \item Player I must ensure that if $Q_n=M$ for some countable $M\prec H((2^\kappa)^+)$, then $p_{n-1},\kappa\in M$.

        \item Player II must ensure that if $Q_n=M$ for some countable $M\prec H((2^\kappa)^+)$, then there exists $\lambda\in\Bf\cap\kappa$ such that $\kappa\cap\hull(M,V_\lambda)\subseteq\lambda$ and $M\proj\lambda\in\M_{p_n}$.

        \item Player I must ensure that if $Q_n=(M,E)$ for some $M\in\M_{p_{n-1}}$ and some $E\in M$, then there exist $M^*\prec (H((2^\kappa)^+),\in,\kappa)$ and a lifting
        $$\pi:M\longrightarrow M^*$$
        such that $p_{n-1}\in\pi(E)$.
        
        \item Player II must ensure that if $Q_n=(M,E)$ for some $M\in\M_{p_{n-1}}$ and some $E\in M$, then there exists $q\in E$ such that $q\in M$ and $p_n\leq q$.
    \end{parts}
    The infinite plays with no rules broken are won by Player II.
    \qed
\end{definition}

We will be interested in those $p$ for which Player II has a winning strategy in $\Game^\mathrm{p}_\kappa(p)$.
Something that we can observe right now, and which will be useful throughout, is that if $\sigma$ is a winning strategy for Player II in $\Game^\mathrm{p}_\kappa(p)$, then any move $q$ that this strategy prescribes for Player II has the same property: Player II has a winning strategy in $\Game^\mathrm{p}_\kappa(q)$.
To see this, one can simply use the tail strategy.

\begin{definition}
    Suppose that $\kappa\in\Bf$.
    We let \intro{$\cc^\mathrm{p}_\kappa$} consist of all $p\in\P^\mathrm{p}\cap V_\kappa$ such that Player II wins $\Game^\mathrm{p}_\kappa(p)$.
    \qed
\end{definition}

The set $\cc^\mathrm{p}_\kappa$ might be empty.
If it is not empty, then it inherits the ordering from $\P^\mathrm{p}$ which makes it into a poset.
It also has the natural unit element.

\begin{notation}
    We denote by \intro{$o$} the pair $(1_\Hi,\emptyset)$.
    \qed
\end{notation}

It turns out that $\cc^\mathrm{p}_\kappa$ is an initial segment of $(\P^\mathrm{p},\leq^\mathrm{p})$.
By this, we simply mean that it is upwards closed.
In particular, if $\cc^\mathrm{p}_\kappa$ contains any element whatsoever, then it contains $o$.
On the other hand, if $\cc^\mathrm{p}_\kappa$ contains $o$, then it also contains $\ran(\sigma)$ for any winning strategy $\sigma$ for Player II in $\Game^\mathrm{p}_\kappa(o)$.

\begin{proposition}\label{258}
    Suppose that $\kappa\in\Bf$.
    Then the following holds.
    \begin{parts}
        \item $\cc^\mathrm{p}_\kappa$ is an initial segment of $(\P^\mathrm{p},\leq^\mathrm{p})$.
        
        \item\label{166} If $\cc^\mathrm{p}_\kappa$ is non-empty, then $o\in\cc^\mathrm{p}_\kappa$ and $(\cc^\mathrm{p}_\kappa,\leq^\mathrm{p},o)$ is a poset.        

        \item\label{168} If $\cc^\mathrm{p}_\kappa$ is non-empty, then, letting $g$ be a $V$-generic for $\cc^\mathrm{p}_\kappa$, we have that the set
        $$\{w_p : p\in g\}$$
        is a $\dd$-generic filter for $\Hi$.
    \end{parts}
\end{proposition}
\begin{proof}
    The fact that $\cc^\mathrm{p}_\kappa$ is an initial segment of $(\P^\mathrm{p},\leq^\mathrm{p})$ easily follows from the fact that we allowed Player II to arbitrarily strengthen condition $p$ at the very beginning of game $\Game^\mathrm{p}_\kappa(p)$.
    Part \ref{166} is then immediate.
    Let us now verify part \ref{168}.
    \begin{pfenum}
        \item Suppose that $\cc^\mathrm{p}_\kappa\not=\emptyset$.

        \item We have that for all $p\in\cc^\mathrm{p}_\kappa$, for all $w\in\Hi$ satisfying that $w\geq w_p$, it holds that $(w,\emptyset)\in\cc^\mathrm{p}_\kappa$ and $(w,\emptyset)\geq p$.
        This readily implies that
        $$\cc^\mathrm{p}_\kappa\Vdash\mbox{``}\{w_p : p\in g\}\mbox{ is a filter in }\Hi\mbox{''}.$$

        \item It remains to verify genericity.
        Let $D\in\dd$ be arbitrary.
        We want to show that
        $$\cc^\mathrm{p}_\kappa\Vdash\exists p\in\dot g, w_p\in D.$$

        \item To this end, it suffices to show that the set
        $$E:=\{p\in\cc^\mathrm{p}_\kappa : \exists w\in D, w_p\leq w\}$$
        is dense in $\cc^\mathrm{p}_\kappa$.

        \item Let $p\in\cc^\mathrm{p}_\kappa$ be arbitrary.
        We want to find $q\in E$ such that $q\leq p$.

        \item Since $p\in\cc^\mathrm{p}_\kappa$, there exists a winning strategy $\sigma$ for Player II in $\Game^\mathrm{p}_\kappa(p)$.

        \item Let us consider the following partial play of $\Game^\mathrm{p}_\kappa(p)$ according to $\sigma$.
        \begin{center}
            \begin{tabular}{c|ccc}
                I  &          & $D$ &      \\
                \hline
                II & $p_{-1}$ &     & $p_0$
            \end{tabular}
        \end{center}

        \item Since $p_{-1}$ and $p_0$ appear in a play according to a winning strategy for Player II, we have that $p_{-1},p_0\in\cc^\mathrm{p}_\kappa$.

        \item The rules of the game imply that
        $$p_0\leq p_{-1}\leq p$$
        and that there exist $w\in D$ such that $w_{p_0}\leq w$.

        \item This imply that $p_0\in E$, so $q:=p_0$ is as required.
    \end{pfenum}
\end{proof}

\begin{proposition}\label{375'}
    Suppose that $\kappa\in\Bf$ and that $\cc^\mathrm{p}_\kappa$ is non-empty.
    Then $\cc^\mathrm{p}_\kappa$ is a proper poset.
\end{proposition}
\begin{prooff}
    \item For $p\in\cc^\mathrm{p}_\kappa$, let $\sigma_p$ be a winning strategy for Player II in the game $\Game^\mathrm{p}_\kappa(p)$.
    
    \item Let 
    $$M\prec (H((2^\kappa)^+),\in,\kappa,\Hi,\dd,(\sigma_p : p\in\cc^\mathrm{p}_\kappa))$$
    be an arbitrary countable virtual model.
    We want to show that $\cc^\mathrm{p}_\kappa$ is proper for $M$.

    \item Let $p\in\cc^\mathrm{p}_\kappa\cap M$ be arbitrary.
    We need to find $p_0\leq p$ which is generic for $(M,\cc^\mathrm{p}_\kappa)$.

    \item Consider the following partial play of $\Game^\mathrm{p}_\kappa(p)$ according to $\sigma_p$.
    \begin{center}
        \begin{tabular}{c|ccc}
            I  &          & $M\proj\kappa$ &       \\
            \hline
            II & $p_{-1}$ &                & $p_0$
        \end{tabular}
    \end{center}
    Since $p_{-1}=\sigma_p(\emptyset)$ and $p\in M$, we get that $p_{-1}\in M\proj\kappa$.
    In other words, Player I did not break the rules.
    
    \item Consequently, it follows that $p_0\leq p_{-1}\leq p$ and that there exists $\lambda\in\Bf\cap\kappa$ such that 
    $$\kappa\cap\hull(P,V_\lambda)\subseteq\lambda$$
    and $M\proj\lambda\in\M_{p_0}$.

    \item We want to show that $p_0$ is as required.
    This comes down to the following.
    Let $q\leq p_0$ and let $E\in M$ be an open dense subset of $\cc^\mathrm{p}_\kappa$ containing $q$.
    We need to show that there exist $r\in E$ such that $r\in M$ and $r\para q$.


    \item Let
    $$\pi:\widehat{M\proj\lambda}\xrightarrow[]{\ \simeq\ }\hull(M,V_\lambda)\prec H_\theta$$
    be the anti-collapse.
    We have that $\pi^{-1}(\kappa)=\lambda$ and $\pi^{-1}(E)=E\cap V_\lambda$.

    \item Let us consider the following partial play of $\Game^\mathrm{p}_\kappa(q)$ according to $\sigma_q$.
    \begin{center}
        \begin{tabular}{c|ccc}
            I  &          & $(M\proj\lambda,E\cap V_\lambda)$ &      \\
            \hline
            II & $q_{-1}$ &         & $q_0$
        \end{tabular}
    \end{center}
    Player I did not break the rules since
    \begin{parts}
        \item $M\proj\lambda\in\M_q\subseteq\M_{q_{-1}}$,

        \item $E\cap V_\lambda=\pi^{-1}(E)\in\pi^{-1}[M]=M\proj\lambda$,

        \item $\pi\rest(M\proj\lambda):M\proj\lambda\to M$ is a lifting,

        \item $q_{-1}\in E=\pi(E\cap V_\lambda)$ ($E$ is open and $q_{-1}\leq q$).
    \end{parts}
    
    \item The rules then imply that $q_0\leq q$ and that there exists $r\in E\cap V_\lambda$ such that $q_0\leq r$ and $r\in M\proj\lambda$.

    \item Thus, $r\in E$, $r\para q$ (as witnessed by $q_0$), and $r\in (M\proj\lambda)\cap V_\lambda\subseteq M$.
\end{prooff}

Thus, if Player II has a winning strategy in $\Game^\mathrm{p}_\kappa(o)$, then $\cc^\mathrm{p}_\kappa$ is a proper poset which adds a $\dd$-generic filter for $\Hi$.
It turns out that a converse is true as well.
To show it, we need the following lemma.

\begin{lemma}
\label{proper consistency}
   The following are equivalent.
    \begin{parts}
        \item\label{proper con a} $(\Hi,\dd)$ is proper- (resp. semiproper-, ssp-) consistent.

        \item\label{proper con b} There exists a proper (resp. semiproper, ssp), complete boolean algebra $\B$ and a mapping $i:\Hi\rightarrow\B$ satisfying that
        \begin{parts}
            \item $i(1)=1$,
            
            \item for all $p,q\in\Hi$, if $p\leq q$, then $i(p)\leq i(q)$,

            \item for all $p,q\in\Hi$, if $p\bot q$, then $i(p)i(q)=0$,

            \item for all $p,q\in\Hi$, if $p\parallel q$, then $i(pq)=i(p)i(q)$,

            \item for all $D\in\dd$, $\sum i[D]=1$.
        \end{parts}
    \end{parts}
\end{lemma}
\begin{proof}
    We explain the proper case, the other two being analogous.
    Let us first show that (\ref{proper con a}$\Rightarrow$\ref{proper con b}).
    Since $(\Hi,\dd)$ is proper-consistent, there exists a proper poset $\Q$ which adds a $\dd$-generic filter for $\Hi$.
    Let $\B$ be the complete boolean algebra obtained as the completion of $\Q$ and let $\dot F$ be a $\B$-name satisfying that
    $$\Vdash_\B\mbox{``}\dot F\mbox{ is }\dd\mbox{-generic for }\Hi\mbox{.''}$$
    Then $i:\Hi\to\B:w\mapsto\|w\in\dot F\|$ is as required.
    To see that (\ref{proper con b}$\Rightarrow$\ref{proper con a}), fix $i:\Hi\to\B$ as in \ref{proper con b}.
    Then for $g$ which is $V$-generic for $\B$, we have that $i^{-1}[g]$ is $\dd$-generic for $\Hi$.
\end{proof}

\begin{theorem}\label{proper con to strategy}
    Suppose that
    \begin{assume}
        \item there exists a proper class of inaccessible cardinals,

        \item $(\Hi,\dd)$ is proper-consistent.
    \end{assume}
    Then for some inaccessible $\kappa$, Player II has a winning strategy in $\Game^\mathrm{p}_\kappa(o)$.
\end{theorem}
\begin{prooff}
    \item By Lemma \ref{proper consistency}, there exists a proper, complete boolean algebra $\B$ and a mapping $i:\Hi\rightarrow\B$ satisfying that
        \begin{parts}
            \item $i(1)=1$,
            
            \item for all $p,q\in\Hi$, if $p\leq q$, then $i(p)\leq i(q)$,

            \item for all $p,q\in\Hi$, if $p\bot q$, then $i(p)i(q)=0$,

            \item for all $p,q\in\Hi$, if $p\parallel q$, then $i(pq)=i(p)i(q)$,

            \item for all $D\in\dd$, $\sum i[D]=1$.
        \end{parts}

    \item Let $\kappa$ be an inaccessible cardinal satisfying that $\rank(\Hi),\rank(\B)<\kappa$ and let us show that Player II has a winning strategy in $\Game^\mathrm{p}_\kappa(o)$.

    \item Let us assume otherwise.
    Since the game is closed for Player II, it follows that Player I has a winning strategy $\sigma$.
    We will reach a contradiction by showing that Player II can defeat $\sigma$.

    \item For $p\in\P^\mathrm{p}$ and $b\in\B$, let $\Psi(p,b)$ be the conjunction of the following statements:
    \begin{parts}[ref=\arabic{pfenumi}$^\circ$\alph{partsi}]
        \item\label{371} for all $M\in\M_{p}$, $\Hi,\B,i\in M$,

        \item\label{373} for all $M\in\M_{p}$, $\lambda_M>\rank(\B)$,

        \item\label{375} $0<b\leq i(w_{p})$,
        
        \item\label{377} for all $M\in\M_{p}$, $b$ is generic for $(M,\B)$.
    \end{parts}
    We will show that Player II can play by maintaining that for all $n\in [-1,\omega)$, there exists $b\in\B$ such that $\Psi(p_n,b)$ holds.

    \item Suppose first that $n=-1$.
    Let 
    $$M\prec (H((2^\kappa)^+),\in,\kappa,\Hi,\B,i)$$
    be countable, let $\lambda:=\sup(\kappa\cap M)$, and let $p_{-1}:=(\emptyset,\{M\proj\lambda\})$.
    We see that conditions \ref{371} and \ref{373} are met.

    \item Since $\B$ is proper and $1_\B\in M$, there exists $b\in\B$ such that
    $$0<b\leq 1=i(w_{p_{-1}})$$
    and such that $b$ is generic for $(M,\B)$.

    \item It is now easily seen that $\Psi(p_{-1},b)$ holds.

    \item Let us now inductively consider the case $n$ assuming that there exists $b\in\B$ such that $\Psi(p_{n-1},b)$ holds.
    Let Player I play a move $Q_n$ and let us show how Player II can answer in such a way as to ensure that there exists $c\in\B$ such that $\Psi(p_{n},c)$ holds.
    

    \item\label{case i}\textbf{Case I.} $Q_n=D$ for some $D\in\dd$.
    \begin{prooff}
        \item Since $\sum i[D]=1$, there exists $w\in D$ such that $i(w)b>0$.

        \item Let $p_n:=(w_{p_{n-1}}w,\M_{p_{n-1}})$ and let $c:=i(w)b$.
        We have that 
        $$0<c\leq i(w_{p_n})$$
        and $\M_{p_n}=\M_{p_{n-1}}$.

        \item Since also $c\leq b$, it follows easily that $\Psi(p_n,c)$ holds.
    \end{prooff}

    \item\textbf{Case II.} $Q_n=M$ for some countable $M\prec H((2^\kappa)^+)$.
    \begin{prooff}
        \item We have that $p_{n-1},\kappa\in M$.

        \item We need to find $p_n\in\P^\mathrm{p}$, $\lambda\in\Bf\cap\kappa$, and $c\in\B$ such that
        \begin{parts}
            \item $p_n\leq p_{n-1}$,

            \item $\kappa\cap\hull(M,V_\lambda)\subseteq\lambda$,
            
            \item $M\proj\lambda\in\M_{p_n}$,
            
            \item $\Psi(p_n,c)$ holds.
        \end{parts}

        \item Let $F:\kappa\to\kappa$ be defined by setting
        $$F(\xi):=\sup(\kappa\cap\hull(M,V_\xi))<\kappa$$
        for $\xi<\kappa$.
        There exists $\lambda\in\Bf\cap\kappa$ such that $F[\lambda]\subseteq\lambda$.

        \item Let $p_n:=(w_{p_{n-1}},\M_{p_{n-1}}\cup\{M\proj\lambda\})\in\P^\mathrm{p}$.
        We have that $p_n\leq p_{n-1}$ and $M\proj\lambda\in\M_{p_n}$.
        It remains to find $c\in\B$ such that $\Psi(p_n,c)$.
        
        \item By elementarity, there exists $b_M\in\B\cap M$ such that $\Psi(p_{n-1},b_M)$ holds.
        
        \item Since $\M_{p_{-1}}\subseteq\M_{p_{n-1}}\in M$, we have that $\B\in M$.
        This means that $\B$ is proper for $M$.

        \item Hence, there exists $c\in\B$ such that $0<c\leq b_M$ and such that $c$ is generic for $(M,\B)$.
        It follows that $\Psi(p_n,c)$ holds, as required.
    \end{prooff}

    \item\textbf{Case III.} $Q_n=(M,E)$ for some $M\in\M_{p_{n-1}}$ and some $E\in M$.
    \begin{prooff}
        \item We have that there exist $M^*\prec(H((2^\kappa)^+),\in,\kappa)$ and a lifting
        $$\pi:M\to M^*$$
        such that $p_{n-1}\in\pi(E)$.

        \item We need to find $p_n\leq p_{n-1}$, $q\in E$, and $c\in\B$ such that
        \begin{parts}
            \item $q\in M$,
            
            \item $p_n\leq q$,

            \item $\Psi(p_n,c)$ holds.
        \end{parts}

        \item Let $\N:=\M_{p_{n-1}}\cap V_{\lambda_M}\in M$ and let $F$ consist of all $d\in\B$ such that there exists $q\in E\cap\P^\mathrm{p}$ satisfying that
        \begin{parts}
            \item for all $P\in\M_q$, $\Hi,\B,i\in P$,

            \item for all $P\in\M_q$, $\lambda_P>\rank(\B)$,

            \item $0<d\leq i(w_q)$,

            \item for all $P\in\M_q$, $d$ is generic for $(P,\B)$,

            \item $\N\subseteq\M_q$.
        \end{parts}
        We have that $F$ is an open subset of $\B-\{0\}$ and $F\in M$.

        \item Let $F^\bot:=\{a\in\B-\{0\} : \forall d\in F,\, ad=0\}$.
        We have that $F\cup F^\bot$ is dense in $\B$ and belongs to $M$.

        \item Since $b$ is generic for $(M,\B)$, there exists $d\in F\cup F^\bot$ such that $bd>0$ and $d\in M$.

        \item Since $p_{n-1}\in\pi(E)$ and $\Psi(p_{n-1},b)$ holds, it follows that $b\in \pi(F)$.

        \item Since $F\subseteq\B$, we have that $F\in M\cap V_{\lambda_M}$.
        This means that $\pi(F)=F$ and consequently, $b\in F$.

        \item Since $bd>0$ and $b\in F$, we have that $d\not\in F^\bot$, which means that $d\in F$.

        \item Since $d\in F\cap M$, there exists $q\in E\cap\P^\mathrm{p}\cap M$ such that
        \begin{parts}
            \item for all $P\in\M_q$, $\Hi,\B,i\in P$,

            \item for all $P\in\M_q$, $\lambda_P>\rank(\B)$,

            \item $0<d\leq i(w_q)$,

            \item for all $P\in\M_q$, $d$ is generic for $(P,\B)$,

            \item $\M_{p_{n-1}}\cap V_{\lambda_M}\subseteq\M_q$.
        \end{parts}

        \item Let $p_n:=(w_{p_{n-1}}w_q,\M_{p_{n-1}}\cup\M_q)$.
        We have that $p_n\in\P^\mathrm{p}$.
        (Note that $0<bd\leq i(w_{p_{n-1}})i(w_q)$, so $w_{p_{n-1}}\para w_q$.)

        \item We also have that $p_n\leq p_{n-1}$ and $p_n\leq q$, while $q\in E\cap M$.

        \item Let $c:=bd>0$.
        We get that $\Psi(p_n,c)$ holds, concluding the argument.
    \end{prooff}

    \item The three cases show together that Player II can survive $\omega$ moves against $\sigma$.
    This contradicts the fact that $\sigma$ is winning for Player I.
\end{prooff}

Hence, we obtain a game-theoretic characterization of proper consistency.

\begin{corollary}\label{878}
    Suppose that there is a proper class of inaccessible cardinals.
    Then the following are equivalent.
    \begin{parts}
        \item $(\Hi,\dd)$ is proper-consistent.

        \item For some inaccessible $\kappa$, Player II has a winning strategy in $\Game^\mathrm{p}_\kappa(o)$.
        \qed
    \end{parts}
\end{corollary}

If $(\Hi,\dd)$ is proper-consistent, then for any sufficiently large inaccessible $\kappa$, the poset $\cc^\mathrm{p}_\kappa$ is a proper poset adding a $\dd$-generic filter for $\Hi$.
It turns out that this is the maximal poset in a collection of posets which are methods identify and which have the the same property.

\begin{proposition}
    Suppose that
    \begin{assume}
        \item $\lambda\in\Bf$ is strictly larger than $\rank(\Hi)$,
        
        \item $\cc$ is a poset contained in $\P^\mathrm{p}\cap V_\lambda$,

        \item for all $p\in\cc$, there exists a winning strategy $\sigma$ for Player II in $\Game^\mathrm{p}_\lambda(p)$ such that $\ran(\sigma)\subseteq\cc$.\footnote{In other words, the moves that $\sigma$ prescribes for Player II are all in $\cc$.}
    \end{assume}
    Then $\cc$ is a proper poset that adds a $\dd$-generic filter for $\Hi$.
    Moreover, $\cc$ is contained in $\cc^\mathrm{p}_\lambda$
\end{proposition}
\begin{proof}
    These assumptions are sufficient to run the proofs of Propositions \ref{258}.\ref{168} and \ref{375'}. 
\end{proof}

There are also minimal ones among the posets that this proposition describes.
They are obtained from by taking the range of a single winning strategy for Player II in $\Game^\mathrm{p}_\lambda(o)$.
We will have the analogous situation in the semiproper and stationary set preserving case.
In Section \ref{An Example for the AS Goodness}, we will actually have to work with a smaller poset in the corresponding collection of posets in the stationary set preserving case.


\section{An Example for the Proper Consistency}
\label{An Example for the Proper Consistency}

In this section, we will illustrate the method described in Section \ref{Proper Consistency} using a concrete example.
First, let us fix our assumptions.

\begin{declaration}
    Suppose that
    \begin{assume}
        \item $X\subseteq\R$ is uncountable,

        \item $\Gamma:[X]^2\to 2$ is an open coloring,

        \item $\kappa$ is an inaccessible,

        \item $<_{V_\kappa}$ is a wellordering of $V_\kappa$.
        \qed
    \end{assume}
\end{declaration}

The goal is to show that there exists a proper poset $\Q$ such that in $V^\Q$, there exists a sequence $(X_n : n<\omega)$ such that $X=\bigsqcup_{n<\omega}X_n$ and such that for all $n<\omega$, $\Gamma\rest [X_n]^2=\mathrm{const}$.\footnote{
This was previously known by the results of \cite{abraham1985on}.
}
This goal can be naturally reformulated as a problem (in our formal sense).

\begin{definition}
    \hfill
    \begin{parts}
        \item Let \intro{$\Hi$} consist of finite partial functions $w\subseteq X\times\Z$ such that for all $n\in\Z$, $\Gamma\rest [w^{-1}(n)]^2=\mathrm{const}$.
        The order on $\Hi$ is the reverse inclusion.

        \item For $x\in X$, let \intro{$D_x$} be $\{w\in\Hi : x\in\dom(w)\}$.

        \item Let \intro{$\dd$} be the set $\{D_x : x\in X\}$.
        \qed
    \end{parts}
\end{definition}

Partial functions $w\in\Hi$ are approximations of the final partition.
We chose to index that partial by $\Z$ in order to simplify notation in arguments below, but of course, any countable set works equally well.
The dense sets inside $\dd$ are such as to ensure that the $\dd$-generic filter for $\Hi$ gives a total partition of $X$.
The following lemma is a simply observation that the problem is what is designed to be.

\begin{lemma}\label{1274}
    Let $\Q$ be an arbitrary poset.
    Then in $V^\Q$, the following are equivalent:
    \begin{parts}
        \item there exists a sequence $(X_n : n<\omega)$ such that $X=\bigsqcup_{n<\omega}X_n$ and such that for all $n<\omega$, $\Gamma\rest [X_n]^2=\mathrm{const}$;

        \item there exists a $\dd$-generic filter for $\Hi$.
        \qed
    \end{parts}
\end{lemma}

This means that our goal reduces to showing that Player II has a winning strategy in $\Game^\mathrm{p}_\kappa(o)$, where the game is defined with respect to the problem $(\Hi,\dd)$.
Our approach for describing this strategy is to define a set of good preconditions in such a way as to be able to show that Player II can win by simply maintaining that all of their moves are good.
The hardest part of this definition is to establish in which way the working part $w_p$ of a good precondition $p$ has to respect its side condition part $\M_p$.
To that end, we first work towards describing a coloring of $X$ ``over'' $M$ given by $\Gamma$, for a countable $M\prec H(\cont^+)$.

\begin{notation}
    We denote by \intro{$<_\hh$} the wellordering $<_{V_\kappa}\rest H(\cont^+)$ of $H(\cont^+)$ and we denote by \intro{$\hh$} the first-order structure
    $$(H(\cont^+),\in,<_\hh, X, \Gamma).$$
    \qed
\end{notation}

The following lemma is a slight generalization of the well-known fact that if $M\prec\hh$ is countable and $\phi(x)$ is a formula with parameters in $M$ and with a realization outside of $M$, then $\phi(x)$ has uncountably many realizations.

\begin{lemma}\label{oca-68}
    Suppose that
    \begin{assume}
        \item $M\prec\hh$ is countable,

        \item $a\subseteq X$ is finite such that $a\cap M=\emptyset$,

        \item $\phi(x)\in\tp^{\hh}(a/M)$.\footnote{
        We use the notation ``$\tp^{\hh}(a/M)$'' to denote the type of the object $a$ in the structure $\hh$ with respect to parameters in $M$.
        }
    \end{assume}
    Then for all countable $S\subseteq X$, there exists a finite $b\subseteq X-S$ such that $\hh\models\phi [b]$.
\end{lemma}
\begin{proof}
    For all countable $S\subseteq X$ satisfying that $S\in M$, it holds that $M\models\phi [a]\wedge a\cap S=\emptyset$.
    Since $M\prec \hh$, it follows that
    $$\hh\models\forall S\in [X]^\omega,\exists x\in [X-S]^{<\omega},\phi(x).$$
    This suffices for the conclusion.
\end{proof}

We can strengthen the above lemma in a way that better takes into the account the coloring $\Gamma$.
Namely, if $a=\{a_i : i<n\}$, we can find fixed colors $(\ge_i : i<n)$ such that for all formulas $\phi$ and all countable $S\subseteq X$ as above, the corresponding tuple $b=\{b_i : i<n\}$ satisfies additionally that the color $\Gamma(a_i,b_i)$ is exactly $\ge_i$ for all $i<n$.

\begin{lemma}\label{oca-86}
    Suppose that
    \begin{assume}
        \item $M\prec\hh$ is countable,

        \item $\sigma\subseteq X-M$ is finite.
    \end{assume}
    Then there exists $\ge:\sigma\to 2$ such that for all $\phi(x)\in\tp^\hh(\sigma/M)$, for all $S\in [X]^\omega$, there exists $f:\sigma\to X-(S\cup\sigma)$ such that
    \begin{parts}
        \item $\hh\models\phi [\ran(f)]$,

        \item for all $a\in \sigma$, $\Gamma(a,f(a))=\ge(a)$.
    \end{parts}
\end{lemma}
\begin{prooff}
    \item Let us assume otherwise.
    Then for all $\ge:\sigma\to 2$, there exist $\phi_\ge(x)\in\tp^\hh(\sigma/M)$ and $S_\ge\in [X]^\omega$ such that for all $f:\sigma\to X-(S_\ge\cup\sigma)$, if $\hh\models\phi_\ge [\ran(f)]$, then for some $a\in\sigma$, $\Gamma(a,f(a))\neq\ge(a)$.
    
    \item\label{oca-104} Let $\phi^\wedge(x):\equiv\bigwedge_{\ge:\sigma\to 2}\phi_\ge(x)$ and let $S:=\bigcup_{\ge:\sigma\to 2}S_\ge\cup\sigma$.
    It follows that for all $f:\sigma\to X-S$, if $\hh\models\phi^\wedge [\ran(f)]$, then
    $$\forall\ge:\sigma\to 2,\exists a\in\sigma,\Gamma(a,f(a))\neq\ge(a).$$

    \item Lemma \ref{oca-68} implies that there exists $f:\sigma\to X-S$ such that $\hh\models\phi^\wedge [\ran(f)]$.

    \item Let $\ge:\sigma\to 2:a\mapsto\Gamma(a,f(a))$.
    The existence of this $\ge$ contradicts \ref{oca-104}.
\end{prooff}

Note that in this lemma, the coloring obtained depends on $\sigma$, so let us call it $\ge_\sigma$.
If we have another finite set $\tau\subseteq X-M$, we would like the colorings $\ge_\sigma$ and $\ge_\tau$ to be coherent, i.e. that $\ge_\sigma\rest (\sigma\cap\tau)=\ge_\tau\rest(\sigma\cap\tau)$.
A priori, there is no reason for this to hold, but the following lemma shows that we can \textit{pick} the colorings in such a way as to ensure it.

\begin{lemma}\label{oca-114}
    There exists a function $M\to\Gamma_M$ which is defined on all countable $M\prec\hh$ and which satisfies that for all countable $M\prec\hh$,
    \begin{parts}
        \item $\Gamma_M:X-M\to 2$,


        \item for all finite $\sigma\subseteq X-M$, for all $\phi(x)\in\tp^\hh(\sigma/M)$, for all $S\in [X]^\omega$, there exists $f:\sigma\to X-(S\cup\sigma)$ such that $\hh\models\phi [\ran(f)]$ and such that for all $a\in \sigma$, $\Gamma(a,f(a))=\Gamma_M(a)$.
    \end{parts}
\end{lemma}
\begin{prooff}
    \item For $M\prec\hh$ countable and for $\sigma\subseteq X-M$ finite, let $\Tilde{\Gamma}(\sigma/M)$ be defined as the $<_\hh$-least $\ge:\sigma\to 2$ such that for all $\phi(x)\in\tp^\hh(\sigma/M)$, for all $S\in [X]^\omega$, there exists $f:\sigma\to X-(S\cup\sigma)$ satisfying that $\hh\models\phi [\ran(f)]$ and that for all $a\in \sigma$, $\Gamma(a,f(a))=\ge(a)$.
    This is possible by Lemma \ref{oca-86}.

    \item Let $\uu$ be an ultrafilter on $[X]^{<\omega}$ satisfying that for all $\sigma\in [X]^{<\omega}$,
    $$\{\tau\in [X]^{<\omega} : \sigma\subseteq\tau\}\in\uu.$$
    
    \item For a countable $M\prec\hh$ and for $\sigma\subseteq X-M$ finite, we define $\Gamma(\sigma/M)$ to be that $\ge:\sigma\to 2$ which satisfies that for $\uu$-almost all $\tau$, it holds that $\sigma\subseteq\tau$ and $\Tilde{\Gamma}(\tau/M)\rest\sigma=\ge$.

    \item Observe that for all countable $M\prec \hh$ and for all finite $\sigma\subseteq\tau\subseteq X-M$, it holds that $\Gamma(\tau/M)\rest\sigma=\Gamma(\sigma/M)$.

    \item For a countable $M\prec\hh$, we let
    $$\Gamma_M:=\bigcup_{\sigma\in [X-M]^{<\omega}}\Gamma(\sigma/M).$$
    It is easily seen that $M\mapsto\Gamma_M$ is as required.
\end{prooff}

\begin{notation}\label{oca-142}
    \hfill
    \begin{parts}
        \item Let us denote by \intro{$M\mapsto\Gamma_M$} the $<_{V_\kappa}$-least mapping whose existence is established by Lemma \ref{oca-114}.

        \item For $M\in\Cor$, we denote by \intro{$\Gamma_M$} the function $\Gamma_{M\cap\hh}$.

        \item We denote by \intro{$\mathrm{r}_2$} the mapping which assigns to each element of $\Z$ its remainder modulo 2.
        \qed
    \end{parts}
\end{notation}

We are now ready to define the set of \textit{good} preconditions.

\begin{definition}
    We define a subset \intro{$\Good$} of $\P^\mathrm{p}\cap V_\kappa$ by asserting that $p\in\Good$ iff
    \begin{parts}
        \item $\M_p\neq\emptyset$,

        \item for all $M\in\M_p$, it holds that $(\Gamma_P : P\prec\hh\mbox{ countable})\in M$,

        \item for all $a\in\dom(w_p)$, it is true that $w_p(a)<0$ if and only if $a$ belongs to the $\in$-least element of $\M_p$,

        \item for all $n\in\Z_{\geq 0}$, for all distinct $a,b\in w_p^{-1}(n)$, it holds that $\Gamma(a,b)=\mathrm{r}_2(n)$,

        \item for all $M,N\in\M_p$, if $M\in N$ and $\M_p\cap (N-M)=\emptyset$, then letting $\sigma:=\dom(w_p)\cap (N-M)$, it holds that $\Gamma_M(\sigma)=(\mathrm{r}_2\circ w_p)\rest\sigma$,

        \item $\dom(w_p)\subseteq\bigcup\M_p$.
        \qed
    \end{parts}
\end{definition}

Recall that the goal is to show that Player II can win $\Game^\mathrm{p}_\kappa(o)$ by maintaining that for all $n\in [-1,\omega)$, $p_n\in\Good$.
The next three lemmas exactly address the three types of moves that Player I can making and show that Player II can keep answering with good preconditions.

\begin{lemma}\label{oca-move one}
    Suppose that $p\in\Good$ and that $a\in X$.
    Then there exists $q\in\Good$ such that $q\leq p$ and $a\in\dom(w_q)$.
\end{lemma}
\begin{prooff}
    \item If $a\in\dom(w_p)$, we can take $q:=p$.
    
    \item Hence, let us assume that $a\not\in \dom(w_p)$.
    
    \item If $a$ in the $\in$-least model of $\M_p$, then let
    \begin{parts}
        \item $n\in\Z_{<0}-\ran(w_p)$,
        
        \item $w_q:=w_p\cup\{(a,n)\}$,
        
        \item $\M_q:=\M_p$,
        
        \item $q:=(w_q,\M_q)$.
    \end{parts}
    It is now easily seen that $q\in\Good$, $q\leq p$, and $a\in\dom(w_q)$.

    \item Hence, let us assume that $a$ is \underline{not} in the $\in$-least model of $\M_p$.

    \item If $a\in\bigcup\M_p$, let
    \begin{parts}
        \item $\M_q:=\M_{p}$,

        \item $M\in N\in\M_p$ be such that $\M_p\cap (N-M)=\emptyset$ and $a\in N-M$,

        \item $n\in\omega-\ran(w_{p})$ be such that $\mathrm{r}_2(n)=\Gamma_M(\{a\})$,

        \item $w_q:=w_{p}\cup\{(a,n)\}$,

        \item $q:=(w_q,\M_q)$.
    \end{parts}
    It is now easily seen that $q\in\Good$, $q\leq p$, and $a\in\dom(w_q)$.

    \item If $a\not\in \bigcup \M_{p}$, let
    \begin{parts}
        \item $N\in\Cor_{<\kappa}$ be such that $\hh,(P\mapsto\Gamma_P),p\in N$,

        \item $\M_q:=\M_{p}\cup\{N\}$,

        \item $M$ be the $\in$-largest element of $\M_p$,
        
        \item $n\in\omega-\ran(w_{p})$ be such that $\mathrm{r}_2(n)=\Gamma_M(\{a\})$,

        \item $w_q:=w_{p}\cup\{(a,n)\}$,

        \item $q:=(w_q,\M_q)$.
    \end{parts}
    It is now easily seen that $q\in\Good$, $q\leq p$, and $a\in\dom(w_q)$.

    \item The four cases conclude the proof.
\end{prooff}

\begin{lemma}\label{oca-move two}
    Suppose that $p\in\Good$ and that $M\prec H((2^\kappa)^+)$ is countable satisfying that $p,\kappa\in M$.
    Then there exist $q\in\Good$ and $\lambda\in\Bf\cap\kappa$ such that $q\leq p$, $\kappa\cap\Hull(M,V_\lambda)\subseteq\lambda$, and $M\proj\lambda\in\M_q$.
\end{lemma}
\begin{proof}
    It is easily seen that we can take $\M_q:=\M_p\cup\{M\proj\lambda\}$ and $q:=(w_p,\M_q)$, for an appropriate $\lambda$.
\end{proof}

\begin{lemma}\label{oca-move three}
    Suppose that
    \begin{assume}
        \item $p\in\Good$,

        \item $N\in\M_p$,

        \item $\pi:N\to N^*\prec H((2^\kappa)^+)$ is a lifting,

        \item $E\in N$ is such that $p\in\pi(E)$.
    \end{assume}
    Then there exist $q\in E\cap N$ and $p'\in\Good$ such that $p'\leq p,q$.
\end{lemma}
\begin{prooff}
    \item Let $n:=|w_p|$ and let $(a_i : i<n)$ be the $<_\hh$-increasing enumeration of $\dom(w_p)$.
    Since $\Gamma$ is an open coloring, there exists a sequence $(t_i : i<n)\in (\omega^{<\omega})^n$ such that
    \begin{parts}
        \item for all $i<n$, $t_i\lhd a_i$,

        \item for all $i<j<n$, $t_i$ and $t_j$ are not initial segments of each other,

        \item for all $i<j<n$, for all $b,c\in X$, if $t_i\lhd b$ and $t_j\lhd c$, then $\Gamma(b,c)=\Gamma(a_i,a_j)$.
    \end{parts}

    \item Let $k\leq n$ be such that $\{a_i : i<k\}=\dom(w_p)\cap N$.

    \item Since $p\in\pi(E)\in N^*$, we may assume without loss of generality that for all $q\in E$, it holds that
    \begin{parts}
        \item $q\in\Good$,
        
        \item $|\dom(w_q)|=n$,

        \item $\M_p\cap N\subseteq \M_q$,

        \item the first $k$ elements of $\dom(w_q)$ are exactly $a_0,\dots,a_{k-1}$,\footnote{
        We mean the first $k$ in the order $<_\hh$.
        }

        \item for all $i<n$, letting $b$ be the $i^\mathrm{th}$ element of $\dom(w_q)$, it holds that $w_q(b)=w_p(a_i)$ and $t_i\lhd b$.
    \end{parts}




    \item\label{oca-305}\claim There exists $q\in E\cap N$ such that for all $i\in [k,n)$, letting $b$ be the $i^\mathrm{th}$ element of $\dom(w_q)$, it holds that $\Gamma(a_i,b)=\mathrm{r}_2(w_p(a_i))$.
    \begin{prooff}
        \item\label{oca-307} Let $\Bar{E}\in N$ be the set of all $n$-tuples that increasingly enumerate $\dom(w_q)$ for some $q\in E$.
        Since $N\prec\widehat{N}$, it suffices to show that there exists $\sigma\in\Bar{E}$ such that for all $i\in [k,n)$, it holds that $\Gamma(a_i,\sigma_i)=\mathrm{r}_2(w_p(a_i))$.
        We may assume without loss of generality that $k<n$, for otherwise there is nothing to show.
        
        \item Let $(M_i : i\leq l)$ be the $\in$-increasing enumeration of $\M_p-N$ and let $(n_i : i\leq l)$ be the increasing sequence of natural numbers such that for all $i\leq l-1$, 
        $$\dom(w_p)\cap (M_{i+1}-M_i)=\{a_j : n_i\leq j<n_{i+1}\}.$$
        Since $k>0$, we must have that $l>0$.

        \item Let $\phi_l(x_k,\dots,x_{n-1})$ be the formula in the language of $\hh$, with parameters in $N\cap H(\cont^+)$, asserting that
        $$a_0^\frown a_1^\frown\cdots^\frown {a_{k-1}}^\frown x_k^\frown \dots^\frown x_{n-1}\in\Bar{E}.$$
        Observe that $\hh\models\phi_l [a_k,\dots,a_{n-1}]$.

        \item By a reverse recursion on $i<l$, we define sequences 
        $$(r_j : n_i\leq j<n_{i+1}),(s_j : n_i\leq j<n_{i+1})\in(\omega^{<\omega})^{n_{i+1}-n_i}$$
        and a formula $\phi_i(x_k,\dots, x_{n_i-1})$.
        We maintain as the recursive hypothesis that 
        $$\hh\models\phi_{i+1}[a_k,\dots,a_{n_{i+1}-1}].$$

        \item Since $p\in\Good$, we have that for all $j\in [n_i,n_{i+1})$, $\Gamma_{M_i}(\{a_j\})=\mathrm{r}_2(w_p(a_j))$.

        \item By the definition of $\Gamma_{M_i}$ and the fact that
        $$\hh\models\phi_{i+1}[a_k,\dots,a_{n_{i+1}-1}],$$
        we have that there exist $b_{n_i},\dots, b_{n_{i+1}-1}\in X-M_{i+1}$ such that
        $$\hh\models\phi_{i+1}[a_k,\dots,a_{n_i-1},b_{n_i},\dots, b_{n_{i+1}-1}]$$
        and such that for all $j\in [n_i,n_{i+1})$, 
        $$\Gamma(a_j, b_j)=\Gamma_{M_i}(\{a_j\})=\mathrm{r}_2(w_p(a_j)).$$

        \item\label{oca-334} For $j\in [n_i,n_{i+1})$, let $r_j\lhd a_j$ and $s_j\lhd b_j$ be finite sequences of natural numbers which are not initial segments of each other and which satisfy that for all $c,d\in X$, if $r_i\lhd c$ and $s_i\lhd d$, then
        $$\Gamma(c,d)=\Gamma(a_j,b_j)=\mathrm{r}_2(w_p(a_j)).$$
        (We are using here that the coloring is open.)

        \item Let $\phi_i(x_k,\dots,x_{n_i-1})$ be the formula asserting that
        $$\exists y_{n_i},\dots,y_{n_{i+1}-1},\bigwedge_{n_i\leq j<n_{i+1}}s_j\lhd y_{j}\wedge$$
        $$\wedge\phi_{i+1}(x_k,\dots,x_{n_i-1},y_{n_i},\dots, y_{n_{i+1}-1}).$$
        This concludes the reverse recursion.

        \item\label{oca-343} Observe that $\hh\models\phi_0$.
        By recursion on $i<l$, we can witness the existential quantifiers in $\phi_i$ so that we obtain $\sigma_{n_i},\dots,\sigma_{n_{i+1}-1}\in X\cap N$ satisfying that $\hh\models\phi_i [\sigma_k,\dots,\sigma_{n_{i+1}-1}]$.

        \item Let $\sigma:=(a_0,\dots,a_{k-1},\sigma_k,\dots,\sigma_{n-1})$.
        Line \ref{oca-343} yields that $\sigma\in\Bar{E}$ and that for all $j\in [k,n)$, $s_j\lhd \sigma_j$.

        \item For all $j\in [k,n)$, the fact that $r_j\lhd a_j$ and $s_j\lhd \sigma_j$ ensures that $\Gamma(a_j,\sigma_j)=\mathrm{r}_2(w_p(a_j))$ (cf. line \ref{oca-334}).
        This shows that $\sigma$ is as required in \ref{oca-307}.
    \end{prooff}

    \item Let $q$ be as provided by Claim \ref{oca-305} and let 
    $$p':=(w_p\cup w_q,\M_p\cup\M_q).$$
    We will show that $q$ and $p'$ are as required.

    \item We already now that $q\in E\cap N$ and that $p'\leq p,q$, so the only thing to verify is that $p'\in\Good$.
    For this, in turn, the only non-trivial point is that for all $n\in\Z_{\geq 0}$, for all distinct $b,c\in w_{p'}^{-1}(n)$, it holds that $\Gamma(b,c)=\mathrm{r}_2(n)$.
    Let us show that this is indeed the case.

    \item This is clear if $b,c\in w_p^{-1}(n)$ or $b,c\in w_q^{-1}(n)$, so let us assume that $b\in w_p^{-1}(n)-N$ and $c\in w_q^{-1}(n)-N$.

    \item Let $i,j\in [k,n)$ be such that $b$ is the $i^\mathrm{th}$ element of $\dom(w_p)$ and $c$ is the $j^\mathrm{th}$ element of $\dom(w_q)$.
    If $i\neq j$, we have that
    $$\Gamma(b,c)=\Gamma(a_i,a_j)=\mathrm{r}_2(n)$$
    because $t_i\lhd b$ and $t_j\lhd c$ (which gives the first equality) and 
    $$w_p(a_i)=w_p(b)=n=w_q(c)=w_p(a_j)$$
    (which gives the seconde equality).

    \item If, on the other hand, $i=j$, then we have that
    $$\Gamma(b,c)=\Gamma(a_i,c)=\mathrm{r}_2(w_p(a_i))=\mathrm{r}_2(n)$$
    because of the way we chose $q$ (cf. Claim \ref{oca-305}).

    \item This shows that $q$ and $p'$ are as required.
\end{prooff}

Putting together the three lemma just proved, we get the desired conclusion.

\begin{corollary}\label{1258}
    Player II has a winning strategy in $\Game^\mathrm{p}_\kappa(o)$.
\end{corollary}
\begin{proof}
    One can easily find $M\in\Cor_{<\kappa}$ such that
    $$p_{-1}:=(\emptyset,\{M\})\in\Good.$$
    Lemmas \ref{oca-move one}, \ref{oca-move two}, and \ref{oca-move three} then show that Player II can survive $\omega$ moves by maintaining that for all $n<\omega$, $p_n\in\Good$.
    This describes a winning strategy for Player II in $\Game^\mathrm{p}_\kappa(o)$.
\end{proof}

As observed in the beginning, this exactly solves the problem that we posed.

\begin{corollary}
    There exists a proper poset $\Q$ such that in $V^\Q$, there exists a sequence
    $$(X_n : n<\omega)$$
    such that $X=\bigsqcup_{n<\omega}X_n$ and for all $n<\omega$, $\Gamma\rest [X_n]^2=\mathrm{const}$.
\end{corollary}
\begin{proof}
    This follows by Corollary \ref{1258} and Lemma \ref{1274}.
\end{proof}


\section[Semiproper Consistency]{Semiproper Consistency\footnote{A version of this section was done jointly with Ben De Bondt. He pursued some of these ideas in a different direction in his thesis \cite{bondt2024stationary}.}}
\label{Semiproper Consistency}

\begin{declaration}
    We fix a problem $(\Hi,\dd)$.
    \qed
\end{declaration}

In this section, we want to characterize the semiproper consistency of $(\Hi,\dd)$ in a way analogous to what we did for the proper consistency.
We will closely follow the outline of Section \ref{Proper Consistency}.
We start by introducing the preconditions.
They are the same as in the proper case, except that we order them differently.

\begin{definition}
    \hfill
    \begin{parts}
        \item The class \intro{$\P^\mathrm{sp}$} consists of all $p=(w_p,\M_p)$ such that $w_p\in\Hi$ and $\M_p$ is a finite vm-chain.

        \item The ordered \intro{$\leq^\mathrm{sp}$} on vm-chains is defined by asserting that $\M\leq^\mathrm{sp}\N$ holds iff for all $N\in\N$, there exists $M\in\M$ such that
        \begin{parts}
            \item $\delta_M=\delta_N$,

            \item $\lambda_M=\lambda_N$,

            \item there exists $X\subseteq V_{\lambda_N}$ such that $M=\Hull(N,X)$.
        \end{parts}

        \item The order \intro{$\leq^\mathrm{sp}$} on $\P^\mathrm{sp}$ is defined by asserting that $p\leq q$ holds iff $w_p\leq_\Hi w_q$ and $\M_p\leq^\mathrm{sp}\M_q$.
        \qed
    \end{parts}
\end{definition}

The game in this case differs from the one in the proper case only in that Player II is required to do less when answering questions of the third kind.

\begin{definition}
    Suppose that $\kappa\in\Bf$ and $p\in\P^\mathrm{sp}\cap V_\kappa$.
    Then the game \intro{$\Game^\mathrm{sp}_\kappa(p)$} is defined as the length $\omega$ two Player game of the form
    \begin{center}
        \begin{tabular}{c|cccccc}
            I  &       & $Q_0$ &       & $Q_1$ &       & $\cdots$\\
            \hline
            II & $p_{-1}$ &       & $p_0$ &       & $p_1$ & $\cdots$
        \end{tabular}
    \end{center}
    where $p_{-1}\in\P^\mathrm{sp}\cap V_\kappa$ is such that $p_{-1}\leq p$, and for all $n<\omega$, the following is satisfied.
    \begin{parts}
        \item Player I must ensure that either
        \begin{parts}
            \item $Q_n=D$ for some $D\in\dd$, or

            \item $Q_n=M$ for some countable $M\prec H((2^\kappa)^+)$, or

            \item $Q_n=(M,E)$ for some $M\in\M_{p_{n-1}}$ and for some $E\in M$.
        \end{parts}

        \item Player II must ensure that $p_n\in\P^\mathrm{sp}\cap V_\kappa$ and $p_n\leq p_{n-1}$.

        \item Player II must ensure that if $Q_n=D$ for some $D\in\dd$, then there exists $w\in D$ such that $w_{p_n}\leq w$.

        \item Player I must ensure that if $Q_n=M$ for some countable $M\prec H((2^\kappa)^+)$, then $p_{n-1},\kappa\in M$.

        \item Player II must ensure that if $Q_n=M$ for some countable $M\prec H((2^\kappa)^+)$, then there exists $\lambda\in\Bf\cap\kappa$ such that $\kappa\cap\hull(M,V_\lambda)\subseteq\lambda$ and $M\proj\lambda\in \M_{p_n}$.

        \item Player I must ensure that if $Q_n=(M,E)$ for some $M\in\M_{p_{n-1}}$ and for some $E\in M$, then there exist $M^*\prec (H((2^\kappa)^+),\in,\kappa)$ and a lifting
        $$\pi:M\longrightarrow M^*$$
        such that $p_{n-1}\in\pi(E)$.
        
        \item Player II must ensure that if $Q_n=(M,E)$ for some $M\in\M_{p_{n-1}}$ and for some $E\in M$, then there exists $q\in E$ such that $\delta(\Hull(M,q))=\delta(M)$ and $p_n\leq q$.
    \end{parts}
    The infinite plays with no rules broken are won by Player II.
    \qed
\end{definition}

\begin{definition}
    Suppose that $\kappa\in\Bf$.
    We let \intro{$\cc^\mathrm{sp}_\kappa$} consist of all $p\in\P^\mathrm{sp}\cap V_\kappa$ such that Player II wins $\Game^\mathrm{sp}_\kappa(p)$.
    \qed
\end{definition}

\begin{proposition}
    Suppose that $\kappa\in\Bf$.
    Then the following holds.
    \begin{parts}
        \item\label{598} The set $\cc^\mathrm{sp}_\kappa$ is an initial segment of $(\P^\mathrm{sp},\leq^\mathrm{sp},o)$.
        
        \item\label{sp166} If $\cc^\mathrm{sp}_\kappa$ is non-empty, then $o\in\cc^\mathrm{sp}_\kappa$ and $(\cc^\mathrm{sp}_\kappa,\leq^\mathrm{sp},o)$ is a poset.        

        \item\label{sp168} If $\cc^\mathrm{sp}_\kappa$ is non-empty, then, letting $g$ be a $V$-generic for $\cc^\mathrm{sp}_\kappa$, we have that the set
        $$\{w_p : p\in g\}$$
        is a $\dd$-generic filter for $\Hi$.
    \end{parts}
\end{proposition}
\begin{proof}
    This is shown analogously to Proposition \ref{258}.
\end{proof}

\begin{proposition}
    Suppose that $\kappa\in\Bf$ and that $\cc^\mathrm{sp}_\kappa$ is non-empty.
    Then $\cc^\mathrm{sp}_\kappa$ is semiproper.
\end{proposition}
\begin{prooff}
    \item For $p\in\cc^\mathrm{sp}_\kappa$, let $\sigma_p$ be a winning strategy for Player II in the game $\Game^\mathrm{sp}_\kappa(p)$.
    
    \item Let
    $$M\prec (H((2^\kappa)^+),\in,\kappa,\Hi,\dd,(\sigma_p : p\in\cc^\mathrm{sp}_\kappa))$$
    be an arbitrary countable virtual model.
    We want to show that $\cc^\mathrm{sp}_\kappa$ is semiproper for $M$.

    \item Let $p\in\cc^\mathrm{sp}_\kappa\cap M$ be arbitrary.
    We need to find $p_0\leq p$ which is semigeneric for $(M,\cc^\mathrm{sp}_\kappa)$.

    \item Consider the following partial play of $\Game^\mathrm{sp}_\kappa(p)$ according to $\sigma_p$.
    \begin{center}
        \begin{tabular}{c|ccc}
            I  &          & $M$ &       \\
            \hline
            II & $p_{-1}$ &                & $p_0$
        \end{tabular}
    \end{center}
    Since $p_{-1}=\sigma_p(\emptyset)$ and $p\in M$, we get that $p_{-1}\in M$.
    Thus, Player I did not break the rules.
    
    \item Consequently, it follows that $p_0\leq p_{-1}\leq p$ and that there exists $\lambda\in\Bf\cap\kappa$ such that 
    $$\kappa\cap\hull(M,V_\lambda)\subseteq\lambda$$
    and such that $M\proj\lambda\in\M_{p_0}$.

    \item We want to show that $p_0$ is as required.
    This comes down to the following.
    Let $q\leq p_0$ and let $E\in M$ be an open dense subset of $\cc^\mathrm{sp}_\kappa$ containing $q$.
    We need to show that there exist $r\in E$ such that $\delta(\Hull(M,r))=\delta(M)$ and $r\para q$.

    \item Let
    $$\pi:\widehat{M\proj\lambda}\xrightarrow[]{\ \simeq\ }\hull(M,V_\lambda)\prec H_\theta$$
    be the anti-collapse.
    We have that $\pi^{-1}(\kappa)=\lambda$ and $\pi^{-1}(E)=E\cap V_\lambda$.

    \item Let us consider the following partial play of $\Game^\mathrm{sp}_\kappa(q)$ according to $\sigma_q$.
    \begin{center}
        \begin{tabular}{c|ccc}
            I  &          & $(N,E\cap V_\lambda)$ &      \\
            \hline
            II & $q_{-1}$ &         & $q_0$
        \end{tabular}
    \end{center}
    Here, $N$ is the unique virtual model in $\M_{q_{-1}}$ such that $\delta_N=\delta_M$.
    (It exists since $M\proj\lambda\in\M_{p_0}$ and $q_{-1}\leq p_0$.)
    
    \item Player I did not break the rules since
    \begin{parts}
        \item $N\in \M_{q_{-1}}$,

        \item $E\cap V_\lambda=\pi^{-1}(E)\in M\proj\lambda\subseteq N$,

        \item $\pi\rest N:N\to\Hull(M,X)\prec H((2^\kappa)^+)$ is a lifting, where $X\subseteq V_\lambda$ is such that $N=\Hull(M\proj\lambda,X)$,

        \item $q_{-1}\in E=\pi(E\cap V_\lambda)$ (because $E$ is open and $q_{-1}\leq q$).
    \end{parts}
    
    \item The rules then imply that $q_0\leq q$ and that there exists $r\in E\cap V_\lambda$ such that $\delta(\Hull(N,r))=\delta(N)$ and $q_0\leq r$.

    \item Thus, $r\in E$, $r\para q$ (as witnessed by $q_0$), and $\delta(\Hull(M,r))=\delta(M)$.
    (For the last point, note that
    $$\delta(M)\leq\delta(\Hull(M,r))=\delta(\Hull(M\proj\lambda,r))\leq$$
    $$\leq\delta(\Hull(N,r))=\delta(N)=\delta(M).)$$
\end{prooff}

We now proceed to showing that the semiproper consistency implies that Player II has a winning strategy in $\Game^\mathrm{sp}_\kappa(o)$, for some $\kappa$.
We will use the following simple lemma.

\begin{lemma}\label{36}
    Suppose that
    \begin{assume}
        \item $\Q$ is a poset,

        \item $\theta\gg\rank(\Q)$ is regular,

        \item $M\prec (H_\theta,\in,\Q)$ is countable,

        \item $p$ is semigeneric for $(M,\Q)$,

        \item $p\in E\in M$.
    \end{assume}
    Then there exist $r\in E\cap\Q$ and $s\leq p,r$ such that $s\Vdash\Check{r}\in\check{M}(\dot{g})$.
    \qed
\end{lemma}

\begin{theorem}\label{semiproper con to strategy}
    Suppose that
    \begin{assume}
        \item there exists a proper class of inaccessible cardinals,

        \item $(\Hi,\dd)$ is semiproper-consistent.
    \end{assume}
    Then for some inaccessible $\kappa$, Player II has a winning strategy in $\Game^\mathrm{sp}_\kappa(o)$.
\end{theorem}
\begin{prooff}
    \item By Lemma \ref{proper consistency}, there exists a semiproper, complete boolean algebra $\B$ and a mapping $i:\Hi\rightarrow\B$ satisfying that
        \begin{parts}
            \item $i(1)=1$,
            
            \item for all $p,q\in\Hi$, if $p\leq q$, then $i(p)\leq i(q)$,

            \item for all $p,q\in\Hi$, if $p\bot q$, then $i(p)i(q)=0$,

            \item for all $p,q\in\Hi$, if $p\parallel q$, then $i(pq)=i(p)i(q)$,

            \item for all $D\in\dd$, $\sum i[D]=1$.
        \end{parts}

    \item Let $\kappa$ be an inaccessible cardinal satisfying that $\rank(\Hi),\rank(\B)<\kappa$ and let us show that $o\in\cc^\mathrm{sp}_\kappa$.
    We need to show that Player II has a winning strategy in $\Game^\mathrm{sp}_\kappa(o)$.

    \item Let us assume otherwise.
    Since the game is closed for Player II, it follows that Player I has a winning strategy $\sigma$.
    We will reach a contradiction by showing that Player II can defeat $\sigma$.

    \item For $p\in\P^\mathrm{sp}$ and $b\in\B$, let $\Psi(p,b)$ be the conjunction of the following statements:
    \begin{parts}[ref=\arabic{pfenumi}$^\circ$\alph{partsi}]
        \item\label{sp371} for all $M\in\M_{p}$, $\Hi,\B,i\in M$,

        \item\label{sp373} for all $M\in\M_{p}$, $\lambda_M>\rank(\B)$,

        \item\label{sp375} $0<b\leq i(w_{p})$,
        
        \item\label{sp377} for all $M\in\M_{p}$, $b$ is semigeneric for $(M,\B)$.
    \end{parts}
    We will show that Player II can play by maintaining that for all $n\in [-1,\omega)$, there exists $b\in\B$ such that $\Psi(p_n,b)$ holds.

    \item Suppose first that $n=-1$.
    Let 
    $$M\prec (H((2^\kappa)^+),\in,\kappa,\Hi,\B,i)$$
    be countable, let $\lambda:=\sup(\kappa\cap M)$, and let $p_{-1}:=(\emptyset,\{M\proj\lambda\})$.
    We see that conditions \ref{sp371} and \ref{sp373} are met.

    \item Since $\B$ is semiproper and $1_\B\in M$, there exists $b\in\B$ such that
    $$0<b\leq 1=i(w_{p_{-1}})$$
    and such that $b$ is semigeneric for $(M,\B)$.

    \item It is now easily seen that $\Psi(p_{-1},b)$ holds.

    \item Let us now inductively consider the case $n$ assuming that there exists $b\in\B$ such that $\Psi(p_{n-1},b)$ holds.
    Let Player I play a move $Q_n$ and let us show how Player II can answer in such a way as to ensure that there exists $c\in\B$ such that $\Psi(p_{n},c)$ holds.
    

    \item\textbf{Case I.} $Q_n=D$ for some $D\in\dd$.
    \begin{proof}
        This is shown like Case I of the proof of Theorem \ref{proper con to strategy}.
    \end{proof}

    \item\textbf{Case II.} $Q_n=M$ for some countable $M\prec H((2^\kappa)^+)$.
    \begin{proof}
        This is shown like Case II of the proof of Theorem \ref{proper con to strategy}.
    \end{proof}

    \item\textbf{Case III.} $Q_n=(M,E)$ for some $M\in\M_{p_{n-1}}$ and some $E\in M$.
    \begin{prooff}
        \item We have that there exist $M^*\prec (H((2^\kappa)^+),\in,\kappa)$ and a lifting
        $$\pi:M\to M^*$$
        such that $p_{n-1}\in\pi(E)$.

        \item\label{860} We need to find $p_n\leq p_{n-1}$, $q\in E$, and $c\in\B$ such that
        \begin{parts}
            \item $\delta(\Hull(M,q))=\delta(M)$,
            
            \item $p_n\leq q$,

            \item $\Psi(p_n,c)$ holds.
        \end{parts}

        \item Let $\N:=\M_{p_{n-1}}\cap V_{\lambda_M}\in M$ and let $F$ consist of all $d\in\B$ such that there exists $q\in E\cap\P^\mathrm{sp}$ satisfying that
        \begin{parts}
            \item for all $N\in\M_q$, $\Hi,\B,i\in N$,

            \item for all $N\in\M_q$, $\lambda_N>\rank(\B)$,

            \item $0<d\leq i(w_q)$,

            \item for all $N\in\M_q$, $d$ is semigeneric for $(N,\B)$,

            \item $\N$ is an initial segment of $(\M_q,\in)$.
        \end{parts}
        We have that $F\in M$.

        \item\claim There exists $d\in F$ and $c\in\B$ such that
        \begin{parts}
            \item for all $N\in\M_{p_{n-1}}-\N$, $\delta(\Hull(N,d))=\delta(N)$,

            \item $0<c\leq bd$,

            \item for all $N\in\M_{p_{n-1}}-\N$, $c$ is semigeneric for $(\Hull(N,d),\B)$.
        \end{parts}
        \begin{prooff}
            \item Since $p_{n-1}\in\pi(E)$ and $\Psi(p_{n-1},b)$ holds, it follows that $b\in \pi(F)$.

            \item Since $F\subseteq\B$, we have that $F\in M\cap V_{\lambda_M}$.
            This means that $\pi(F)=F$ and consequently, $b\in F$.

            \item By Lemma \ref{36}, there exists $d\in F$ and $c\in\B$ such that
            $$0<c\leq bd$$
            and $c\Vdash d\in M(\dot g)$.

            \item Let $N\in\M_{p_{n-1}}-\N$ be arbitrary.
            We have to verify that $\delta(\Hull(N,d))=\delta(N)$ and $c$ is semigeneric for $(\Hull(N,d),\B)$.

            \item\label{904} Since $c\leq b$ and $b$ is semigeneric for $N$, we have that
            $$c\Vdash\delta(N(\dot g))=\delta(N).$$

            \item Since $c\Vdash d\in M(\dot g)$, we have that $c\Vdash d\in N(\dot g)$
            
            \item\label{909} Consequently, $c\Vdash N\prec \Hull(N,d)\prec N(\dot g)$.

            \item Lines \ref{904} and \ref{909} imply that $\delta(\Hull(N,d))=\delta(N)$.


            \item Line \ref{909} also implies that $c\Vdash\Hull(N,d)(\dot g)=N(\dot g)$.
            
            \item By referencing line \ref{904} once again, we get that $c$ is semigeneric for $(\Hull(N,d),\B)$.
        \end{prooff}

        \item Since $d\in F\cap\Hull(M,d)$, there exists $q\in E\cap\P^\mathrm{sp}\cap\Hull(M,d)$ such that
        \begin{parts}
            \item for all $N\in\M_q$, $\Hi,\B,i\in N$,

            \item for all $N\in\M_q$, $\lambda_N>\rank(\B)$,

            \item $0<d\leq i(w_q)$,

            \item for all $N\in\M_q$, $d$ is semigeneric for $(N,\B)$,

            \item $\M_{p_{n-1}}\cap V_{\lambda_M}$ is an initial segment of $(\M_q,\in)$.
        \end{parts}

        \item Let
        \begin{parts}
            \item $w_{p_n}:=w_{p_{n-1}}w_q$,

            \item $\M_{p_n}:=\M_q\cup\{\Hull(N,d) : N\in\M_{p_{n-1}},\, \delta_N\geq\delta_M\}$,
            
            \item $p_n:=(w_{p_n},\M_{p_n})$.
        \end{parts}

        \item \claim $p_n$, $q$, and $c$ are as required in \ref{860}
        \begin{prooff}
            \item Since $0<c\leq bd\leq i(w_{p_{n-1}})i(w_q)$, it follows that $w_{p_{n-1}}\para w_q$.
            This means that $w_{p_n}$ is a well defined element of $\Hi$.

            \item Since $\M_q\in\Hull(M,d)$, we can easily verify that $\M_{p_n}$ is a vm-chain.

            \item This means that $p_n\in\P^\mathrm{sp}$.

            \item The fact that $p_n\leq p_{n-1}$ is self-evident, modulo the observation that for all $N\in\M_{p_{n-1}}-\N$, we have that
            $$\delta(\Hull(N,d))=\delta(N)\mbox{ and } \lambda(\Hull(N,d))=\lambda(N).$$

            \item $q$ was picked so that $q\in E$, while the fact that $q\in\Hull(M,d)$ easily implies that
            $$\delta(\Hull(M,q))=\delta(M).$$

            \item The fact that $p_n\leq q$ is also self-evident.

            \item It remains to verify that $\Psi(p_n,c)$ holds, i.e. that
            \begin{parts}
                \item\label{958} for all $N\in\M_{p_n}$, $\Hi,\B,i\in N$,

                \item\label{960} for all $N\in\M_{p_n}$, $\lambda_N>\rank(\B)$,
        
                \item\label{962} $0<c\leq i(w_{p_n})$,
                
                \item\label{964} for all $N\in\M_{p_n}$, $c$ is semigeneric for $(N,\B)$.
            \end{parts}
            Parts \ref{958} and \ref{960} are immediate, while for \ref{962}, we have
            $$c\leq bd\leq i(w_{p_{n-1}})i(w_q)=i(w_{p_n}).$$

            \item Let us explain part \ref{964}.
            By the choice of $q$, we have that $d$ is semigeneric for $(N,\B)$ for all $N\in\M_q$.

            \item By the choice of $c$, we know that it is semigeneric for 
            $$(\Hull(N,d),\B),$$
            for all $N\in\M_q-\N$.

            \item Since $c\leq d$, the last two line imply that $c$ is generic for $(N,\B)$ for all $N\in\M_{p_n}$.
        \end{prooff}
        
        \item The claim concludes verification of Case III.
    \end{prooff}

    \item The three cases show together that Player II can survive $\omega$ moves against $\sigma$.
    This contradicts the fact that $\sigma$ is winning for Player I.
\end{prooff}

\begin{corollary}\label{1702}
    Suppose that there is a proper class of inaccessible cardinals.
    Then the following are equivalent.
    \begin{parts}
        \item $(\Hi,\dd)$ is semiproper-consistent.

        \item For some inaccessible $\kappa$, Player II has a winning strategy in $\Game^\mathrm{sp}_\kappa(o)$.
        \qed
    \end{parts}
\end{corollary}


\section{Stationary Set Preserving Consistency}
\label{ssp consistency}

\begin{declaration}
    We fix a problem $(\Hi,\dd)$.
    \qed
\end{declaration}

In this section, we want to characterize the ssp consistency of $(\Hi,\dd)$ in a way analogous to what we did for proper and semiproper consistency.
We follow closely the outline of Sections \ref{Proper Consistency} and \ref{Semiproper Consistency}.
In this case, the preconditions and their ordering are exactly the same as in the semiproper case.
The game $\Game^\mathrm{ssp}_\kappa(p)$ that we define below differs from the game $\Game^\mathrm{sp}_\kappa(p)$ in the questions of the second type that Player I can ask.
While a semiproper poset is semiproper for club many models, this is not true for an ssp poset.
What we have is that an ssp poset is semiproper for projectively stationary many models (cf. \cite[Lemma 4.8]{feng2003structure}).
The change in the game is such that it reflects this fact.

\begin{definition}
    Suppose that $\kappa\in\Bf$ and $p\in\P^\mathrm{sp}\cap V_\kappa$.
    Then the game \intro{$\Game^\mathrm{ssp}_\kappa(p)$} is defined as the length $\omega$ two Player game of the form
    \begin{center}
        \begin{tabular}{c|cccccc}
            I  &       & $Q_0$ &       & $Q_1$ &       & $\cdots$\\
            \hline
            II & $p_{-1}$ &       & $p_0$ &       & $p_1$ & $\cdots$
        \end{tabular}
    \end{center}
    where $p_{-1}\in\P^\mathrm{sp}\cap V_\kappa$ is such that $p_{-1}\leq p$ and for all $n<\omega$, the following is satisfied.
    \begin{parts}
        \item Player I must ensure that either
        \begin{parts}
            \item $Q_n=D$ for some $D\in\dd$, or

            \item $Q_n=(U,S)$ for some $U\subseteq H((2^\kappa)^+)$ and for some $S\subseteq\omega_1$ which is stationary,

            \item $Q_n=(M,E)$ for some $M\in\M_{p_{n-1}}$ and for some $E\in M$.
        \end{parts}

        \item Player II must ensure that $p_n\in\P^\mathrm{sp}\cap V_\kappa$ and $p_n\leq p_{n-1}$.

        \item Player II must ensure that if $Q_n=D$ for some $D\in\dd$, then there exists $w\in D$ such that $w_{p_n}\leq w$.

        \item Player II must ensure that if $Q_n=(U,S)$ for some $U\subseteq H((2^\kappa)^+)$ and for some $S\subseteq\omega_1$ which is stationary, then there exist
        $$M\prec (H((2^\kappa)^+),\in,\kappa,p_{n-1},U)$$
        and $\lambda\in\Bf\cap\kappa$ such that 
        $$\kappa\cap\hull(M,V_\lambda)\subseteq\lambda,$$
        $M\proj\lambda\in \M_{p_n}$, and $\delta_M\in S$.

        \item Player I must ensure that if $Q_n=(M,E)$ for some $M\in\M_{p_{n-1}}$ and for some $E\in M$, then there exist $M^*\prec (H((2^\kappa)^+),\in,\kappa)$ and a lifting
        $$\pi:M\longrightarrow M^*$$
        such that $p_{n-1}\in\pi(E)$.
        
        \item Player II must ensure that if $Q_n=(M,E)$ for some $M\in\M_{p_{n-1}}$ and for some $E\in M$, then there exists $q\in E$ such that $\delta(\Hull(M,q))=\delta(M)$ and $p_n\leq q$.
    \end{parts}
    The infinite plays with no rules broken are won by Player II.
    \qed
\end{definition}

\begin{definition}
    Suppose that $\kappa\in\Bf$.
    We let \intro{$\cc^\mathrm{ssp}_\kappa$} consist of all $p\in\P^\mathrm{sp}\cap V_\kappa$ such that Player II wins $\Game^\mathrm{ssp}_\kappa(p)$.
    \qed
\end{definition}

\begin{proposition}\label{1129}
    Suppose that $\kappa\in\Bf$.
    Then the following holds.
    \begin{parts}
        \item\label{s598} The set $\cc^\mathrm{ssp}_\kappa$ is an initial segment of $(\P^\mathrm{sp},\leq^\mathrm{sp},o)$.
        
        \item\label{ssp166} If $\cc^\mathrm{ssp}_\kappa$ is non-empty, then $o\in\cc^\mathrm{ssp}_\kappa$ and $(\cc^\mathrm{ssp}_\kappa,\leq^\mathrm{sp},o)$ is a poset.        

        \item\label{ssp168} If $\cc^\mathrm{ssp}_\kappa$ is non-empty, then, letting $g$ be a $V$-generic for $\cc^\mathrm{ssp}_\kappa$, we have that the set
        $$\{w_p : p\in g\}$$
        is a $\dd$-generic filter for $\Hi$.
    \end{parts}
\end{proposition}
\begin{proof}
    This is verified as in the proof of Proposition \ref{258}.
\end{proof}

\begin{proposition}\label{1146}
    Suppose that $\kappa\in\Bf$ and that $\cc^\mathrm{ssp}_\kappa$ is non-empty.
    Then $\cc^\mathrm{ssp}_\kappa$ is stationary set preserving.
\end{proposition}
\begin{prooff}
    \item\label{1073} Let us assume otherwise.
    Then there exist $S\subseteq\omega_1$ which is stationary, $p\in\ssp$, and $\dot{C}$ which is a canonical name for a subset of $\omega_1$ such that $p\Vdash$``$\Dot{C}$ is a club and $\check{S}\cap\dot{C}=\emptyset$''.
    
    \item For all $q\in\cc^\mathrm{ssp}_\kappa$, let $\sigma_q$ be a winning strategy for Player II in the game $\Game^\mathrm{ssp}_\kappa(p)$.

    \item Consider the following partial play of $\Game^\mathrm{ssp}_\kappa(p)$ according to $\sigma_p$.
    \begin{center}
        \begin{tabular}{c|ccc}
            I  &          & $(U,S)$ &       \\
            \hline
            II & $p_{-1}$ &                & $p_0$
        \end{tabular}
    \end{center}
    Here, $U\subseteq H((2^\kappa)^+)$ is some canonical coding of the tuple
    $$(\Hi,\dd,(\sigma_p : p\in\cc^\mathrm{ssp}_\kappa),\dot{C}).$$

    \item The rules then imply that $p_0\leq p$ and that for some countable
    $$M\prec (H((2^\kappa)^+),\in,\kappa,p_{-1},\Hi,\dd,(\sigma_p : p\in\cc^\mathrm{ssp}_\kappa),\dot{C})$$
    and for some $\lambda\in\Bf\cap\kappa$, we have that
    $$\kappa\cap\hull(M,V_\lambda)\subseteq\lambda,$$
    $M\proj\lambda\in \M_{p_0}$, and $\delta_M\in S$.

    \item\label{1173} \claim $p_0$ is semigeneric for $(M,\ssp)$.
    \begin{prooff}
        \item We need to verify the following.
        Let $q\in\cc^\mathrm{ssp}_\kappa$ be such that $q\leq p_0$ and let $E$ be a dense open subset of $\cc^\mathrm{ssp}_\kappa$ such that $q\in E\in M$.
        We need to find $r\in E$ such that $\delta(\Hull(M,r))=\delta(M)$ and $r\para q$.

        \item Let us consider the following partial play of $\Game^\mathrm{ssp}_\kappa(q)$ according to $\sigma_q$.
        \begin{center}
            \begin{tabular}{c|ccc}
                I  &          & $(N,E\cap V_\lambda)$ &      \\
                \hline
                II & $q_{-1}$ &         & $q_0$
            \end{tabular}
        \end{center}
        Here, $N$ is the unique virtual model in $\M_{q_{-1}}$ such that $\delta_N=\delta_M$.
        (It exists since $M\proj\lambda\in\M_{p_0}$ and $q_{-1}\leq p_0$.)

        \item Let
        $$\pi:\widehat{M\proj\lambda}\xrightarrow[]{\ \simeq\ }\hull(M,V_\lambda)\prec H_\theta$$
        be the anti-collapse.
        We have that $\pi^{-1}(\kappa)=\lambda$ and $\pi^{-1}(E)=E\cap V_\lambda$.
        
        \item Player I did not break the rules since
        \begin{parts}
            \item $N\in \M_{q_{-1}}$,
    
            \item $E\cap V_\lambda=\pi^{-1}(E)\in M\proj\lambda\subseteq N$,
    
            \item $\pi\rest N:N\to\Hull(M,X)\prec H((2^\kappa)^+)$ is a lifting, where $X\subseteq V_\lambda$ is such that $N=\Hull(M\proj\lambda,X)$,
    
            \item $q_{-1}\in E=\pi(E\cap V_\lambda)$ (because $E$ is open and $q_{-1}\leq q$).
        \end{parts}
        
        \item The rules then imply that $q_0\leq q$ and that there exists $r\in E\cap V_\lambda$ such that $\delta(\Hull(N,r))=\delta(N)$ and $q_0\leq r$.
    
        \item Thus, $r\in E$, $r\para q$ (as witnessed by $q_0$), and $\delta(\Hull(M,r))=\delta(M)$.
        (For the last point, note that
        $$\delta(M)\leq\delta(\Hull(M,r))=\delta(\Hull(M\proj\lambda,r))\leq$$
        $$\leq\delta(\Hull(N,r))=\delta(N)=\delta(M).)$$
    \end{prooff}

    \item Now, let $g$ be a $V$-generic filter for $\ssp$ that contains $p_0$.
    Since $p_0$ is semigeneric for $(M,\ssp)$, we have that
    $$\delta_{M[g]}=\delta_M.$$

    \item Since $\dot{C}\in M$ and $\dot{C}_g$ is a club in $\omega_1^V$, we have that
    $$\delta_M=\delta_{M[g]}\in\dot C_g.$$

    \item Since $\delta_M\in S$, we get that $\dot C_g\cap S\not=\emptyset$, which contradicts line \ref{1073} and the fact that $p_0\leq p$.
\end{prooff}

We now turn towards showing that the ssp consistency implies the existence of a winning strategy for Player II in $\Game^\mathrm{ssp}_\kappa(o)$, for some $\kappa$.
We will need the following lemma.

\begin{lemma}\label{ssp and sp}
    Suppose that
    \begin{assume}
        \item $\Q$ is a poset,

        \item $\theta\gg\rank(\Q)$ is regular,

        \item $\Q$ is stationary set preserving,

        \item $S\subseteq\omega_1$ is stationary,

        \item $U$ is an arbitrary subset of $H_\theta$.
    \end{assume}
    Then there exists a countable $M\prec (H_\theta,\in,U)$ such that $\delta_M\in S$ and $\Q$ is semiproper for $M$.
\end{lemma}
\begin{proof}
    See \cite[Lemma 4.8]{feng2003structure}.
\end{proof}

\begin{theorem}
    Suppose that
    \begin{assume}
        \item there exists a proper class of inaccessible cardinals,

        \item $(\Hi,\dd)$ is ssp-consistent.
    \end{assume}
    Then for some inaccessible $\kappa$, Player II has a winning strategy in $\Game^\mathrm{ssp}_\kappa(o)$.
\end{theorem}
\begin{prooff}
    \item By Lemma \ref{proper consistency}, there exists an ssp, complete boolean algebra $\B$ and a mapping $i:\Hi\rightarrow\B$ satisfying that
        \begin{parts}
            \item for all $w\in\Hi$, $i(w)>0$,
            
            \item $i(1)=1$,
            
            \item for all $p,q\in\Hi$, if $p\leq q$, then $i(p)\leq i(q)$,

            \item for all $p,q\in\Hi$, if $p\bot q$, then $i(p)i(q)=0$,

            \item for all $p,q\in\Hi$, if $p\parallel q$, then $i(pq)=i(p)i(q)$,

            \item for all $D\in\dd$, $\sum i[D]=1$.
        \end{parts}

    \item Let $\kappa$ be an inaccessible cardinal satisfying that $\rank(\Hi),\rank(\B)<\kappa$ and let us show that $o\in\cc^\mathrm{ssp}_\kappa$.
    We need to show that Player II has a winning strategy in $\Game^\mathrm{ssp}_\kappa(o)$.

    \item Let us assume otherwise.
    Since the game is closed for Player II, it follows that Player I has a winning strategy $\sigma$.
    We will reach a contradiction by showing that Player II can defeat $\sigma$.

    \item For $p\in\P^\mathrm{sp}$ and $b\in\B$, let $\Psi(p,b)$ be the conjunction of the following statements:
    \begin{parts}[ref=\arabic{pfenumi}$^\circ$\alph{partsi}]
        \item\label{ssp371} for all $M\in\M_{p}$, $\Hi,\B,i\in M$,

        \item\label{ssp373} for all $M\in\M_{p}$, $\lambda_M>\rank(\B)$,

        \item\label{ssp375} $0<b\leq i(w_{p})$,
        
        \item\label{ssp377} for all $M\in\M_{p}$, $b$ is semigeneric for $(M,\B)$.
    \end{parts}
    We will show that Player II can play by maintaining that for all $n\in [-1,\omega)$, there exists $b\in\B$ such that $\Psi(p_n,b)$ holds.

    \item Suppose first that $n=-1$.
    By Lemma \ref{ssp and sp}, there exists a countable
    $$M\prec (H((2^\kappa)^+),\in,\kappa,\Hi,\B,i)$$
    such that $\B$ is semiproper for $M$.

    \item Let $\lambda:=\sup(\kappa\cap M)$ and let $p_{-1}:=(\emptyset,\{M\proj\lambda\})$.
    We see that conditions \ref{ssp371} and \ref{ssp373} are met.

    \item Since $\B$ is semiproper for $M$ and $1_\B\in M$, there exists $b\in\B$ such that
    $$0<b\leq 1=i(w_{p_{-1}})$$
    and such that $b$ is semigeneric for $(M,\B)$.
    It is easily seen that $\Psi(p_{-1},b)$ holds.

    \item Let us now inductively consider the case $n$ assuming that there exists $b\in\B$ such that $\Psi(p_{n-1},b)$ holds.
    Let Player I make a move $Q_n$ and let us show how Player II can answer in such a way as to ensure that there exists $c\in\B$ such that $\Psi(p_{n},c)$ holds.
    

    \item\textbf{Case I.} $Q_n=D$ for some $D\in\dd$.
    \begin{proof}
        This is verified in the same way as Case I of the proof of Theorem \ref{proper con to strategy}.
    \end{proof}

    \item\textbf{Case II.} $Q_n=(U,S)$ for some $U\subseteq H((2^\kappa)^+)$ and for some $S\subseteq\omega_1$ which is stationary.
    \begin{prooff}
        \item By Lemma \ref{ssp and sp}, there exists
        $$M\prec(H((2^\kappa)^+),\in,\kappa,p_{n-1},U)$$
        such that $\delta_M\in S$ and $\B$ is semiproper for $M$.

        \item Let $F:\kappa\to\kappa$ be defined by setting
        $$F(\xi):=\sup(\kappa\cap\hull(M,V_\xi))<\kappa$$
        for $\xi<\kappa$.
        There exists $\lambda\in\Bf\cap\kappa$ such that $F[\lambda]\subseteq\lambda$.

        \item Let $p_n:=(w_{p_{n-1}},\M_{p_{n-1}}\cup\{M\proj\lambda\})\in\P^\mathrm{sp}$.
        We have that $p_n\leq p_{n-1}$ and $M\proj\lambda\in\M_{p_n}$.
        It remains to find $c\in\B$ such that $\Psi(p_n,c)$.

        \item By elementarity, there exists $b_M\in\B\cap M$ such that $\Psi(p_{n-1},b_M)$ holds.
        
        \item Since $\B$ is semiproper for $M$, there exists $c\in\B$ such that $0<c\leq b_M$ and such that $c$ is semigeneric for $(M,\B)$.
        It follows that $\Psi(p_n,c)$ holds, as required.
    \end{prooff}

    \item\textbf{Case III.} $Q_n=(M,E)$ for some $M\in\M_{p_{n-1}}$ and some $E\in M$.
    \begin{proof}
        This is shown like Case III of the proof of Theorem \ref{semiproper con to strategy}.
    \end{proof}

    \item The three cases show together that Player II can survive $\omega$ moves against $\sigma$.
    This contradicts the fact that $\sigma$ is winning for Player I.
\end{prooff}

\begin{corollary}\label{1993}
    Suppose that there is a proper class of inaccessible cardinals.
    Then the following are equivalent.
    \begin{parts}
        \item $(\Hi,\dd)$ is ssp-consistent.

        \item For some inaccessible $\kappa$, Player II has a winning strategy in $\Game^\mathrm{ssp}_\kappa(o)$.
        \qed
    \end{parts}
\end{corollary}

\begin{proposition}\label{1806}
    Suppose that
    \begin{assume}
        \item $\lambda\in\Bf$ is strictly larger than $\rank(\Hi)$,
        
        \item $\cc$ is a poset contained in $\P^\mathrm{sp}\cap V_\lambda$,

        \item for all $p\in\cc$, there exists a winning strategy $\sigma$ for Player II in $\Game^\mathrm{ssp}_\lambda(p)$ such that $\ran(\sigma)\subseteq\cc$.
    \end{assume}
    Then $\cc$ is a stationary set preserving poset that adds a $\dd$-generic filter for $\Hi$.
\end{proposition}
\begin{proof}
    These assumptions are sufficient to run the proofs of Propositions \ref{1129}.\ref{ssp168} and \ref{1146}.
\end{proof}


\section{AS Goodness}
\label{AS Goodness}

\begin{declaration}
    We fix a problem $(\Hi,\dd)$.
    \qed
\end{declaration}

We want to give a sufficient condition for this problem to be ssp-consistent.
This is the content of the following definition.

\begin{definition}
    Suppose that $\kappa>\rank(\Hi)$ is inaccessible.
    Then $(\Hi,\dd)$ is \intro{AS-good} at $\kappa$ iff for all stationary $S\subseteq\omega_1$, it holds in $V^{\Col(\omega,<\kappa)}$ that for all $\dd$-generic filters $F$ for $\Hi$, there exists an elementary embedding
    $$j:V\to M$$
    such that $\crit(j)=\omega_1^V\in j(S)$ and there exists a $j(\dd)$-generic filter $F^*$ for $j(\Hi)$ such that $j[F]\subseteq F^*$.
    \qed
\end{definition}

In Section \ref{An Example for the AS Goodness}, we will illustrate how this criterion can be used in practice.
We now proceed to show that it implies the ssp consistency.
More precisely, we will show how Player II can win $\Game^\mathrm{ssp}_\kappa(o)$.
In the next definition, we recursively define a family of posets
$$(\cc^\mathrm{BV}_\lambda : \lambda\in\Bf,\, \lambda>\rank(\Hi))$$
which is closely related to the posets $\cc^\mathrm{ssp}_\lambda$.
The conditions in these posets are chosen based on a modification of the game $\Game^\mathrm{ssp}_\kappa(p)$: Player I is not allowed any more to play the moves $Q_n=(U,S)$, while in the moves $Q_n=(M,E)$, he needs to meet a slightly different requirement.

\begin{definition}\label{1488}
    Suppose that
    \begin{assume}
        \item $\eta\in\Bf$ is strictly larger than $\rank(\Hi)$,

        \item for all $\xi\in\Bf\cap\eta$, the poset $\cc^\mathrm{BV}_\xi$ has been defined,

        \item for all $\xi\in\Bf\cap\eta$, poset $\cc^\mathrm{BV}_\xi$ is an initial segment of $\P^\mathrm{sp}\cap V_\xi$.
    \end{assume}
    For $p\in\P^\mathrm{sp}\cap V_\eta$, the game \intro{$\Game^\mathrm{BV}_\eta(p)$} is defined as the length $\omega$ two Player game of the form
    \begin{center}
        \begin{tabular}{c|cccccc}
            I  &       & $Q_0$ &       & $Q_1$ &       & $\cdots$\\
            \hline
            II & $p_{-1}$ &       & $p_0$ &       & $p_1$ & $\cdots$
        \end{tabular}
    \end{center}
    where $p_{-1}\in\P^\mathrm{sp}\cap V_\eta$ is such that $p_{-1}\leq p$ and for all $n<\omega$, the following is satisfied.
    \begin{parts}
        \item Player I must ensure that either
        \begin{parts}
            \item $Q_n=D$ for some $D\in\dd$, or


            \item $Q_n=(M,E)$ for some $M\in\M_{p_{n-1}}$ and for some $E\in M$.
        \end{parts}

        \item Player II must ensure that $p_n\in\P^\mathrm{sp}\cap V_\eta$ and $p_n\leq p_{n-1}$.

        \item Player II must ensure that if $Q_n=D$ for some $D\in\dd$, then there exists $w\in D$ such that $w_{p_n}\leq w$.


        \item Player I must ensure that if $Q_n=(M,E)$ for some $M\in\M_{p_{n-1}}$ and for some $E\in M$, then $E$ is a dense subset $\cc^\mathrm{BV}_{\lambda_M}$.
        
        \item Player II must ensure that if $Q_n=(M,E)$ for some $M\in\M_{p_{n-1}}$ and for some $E\in M$, then there exists $q\in E$ such that $\delta(\Hull(M,q))=\delta(M)$ and $p_n\leq q$.
    \end{parts}
    The infinite plays with no rules broken are won by Player II.
    The set \intro{$\cc^\mathrm{BV}_\eta$} is defined to consist of all $p\in\P^\mathrm{sp}\cap V_\eta$ such that Player II has a winning strategy in $\Game^{\mathrm{BV}}_\eta(p)$.\footnote{
    This definition is essentially due to De Bondt and Veli\v ckovi\'c.
    }
    \qed
\end{definition}

Note that this game is local: there is no essential difference between $\Game^\mathrm{BV}_\eta(p)$ and $\Game^\mathrm{BV}_{\eta'}(p)$, as long as $\eta,\eta'>\rank(p)$.
We will show that Player II can win $\Game^\mathrm{ssp}_\kappa(o)$ by maintaining that for each of his moves $p_n$, he can win the (auxiliary) game $\Game^\mathrm{BV}_\kappa(p_n)$.
The following lemma observes some basic facts about the posets $\cc^\mathrm{BV}_\eta$, which are standard by now.

\begin{lemma}\label{1535}
    Suppose that $\eta\in\Bf$ is strictly larger than $\rank(\Hi)$.
    The following holds.
    \begin{parts}
        \item The poset $\cc^\mathrm{BV}_\eta$ is well-defined.
        
        \item For all $w\in\Hi$, $(w,\emptyset)\in\cc^\mathrm{BV}_\eta$.

        \item The poset $\cc^\mathrm{BV}_\eta$ is an initial segment of $\P^\mathrm{sp}\cap V_\eta$.

        \item For all $\zeta\in\Bf$ which are larger than $\eta$, $\cc^\mathrm{BV}_\zeta\cap V_\eta=\cc^\mathrm{BV}_\eta$.

        \item If $g$ is $V$-generic for $\cc^\mathrm{BV}_\eta$, then $\{w_p : p\in g\}$ is a $\dd$-generic filter for $\Hi$.
        \qed
    \end{parts}
\end{lemma}

As we said, Player II will win $\Game^\mathrm{ssp}_\kappa(o)$ by maintaining that he can win $\Game^\mathrm{BV}_\kappa(p_n)$ for all $n<\omega$.
The following proposition shows how he can answer a question of the form $Q_n=(U,S)$ of Player I in $\Game^\mathrm{ssp}_\kappa(o)$.
The proof uses the AS goodness in a crucial way.

\begin{proposition}\label{134}
    Suppose that
    \begin{assume}
        \item $\kappa>\rank(\Hi)$ is inaccessible,

        \item $(\Hi,\dd)$ is AS-good at $\kappa$,

        \item $p\in\cc^\mathrm{BV}_\kappa$,

        \item $U\subseteq H((2^\kappa)^+)$,

        \item $S\subseteq\omega_1$ is stationary.
    \end{assume}
    Then there exists a countable 
    $$M\prec (H((2^\kappa)^+),\in,U,\kappa,p)$$
    with $\delta_M\in S$, there exits $\lambda\in\Bf\cap\kappa$ with $\Hull(M,V_\lambda)\cap\kappa\subseteq\lambda$, and there exists $q\in\cc^\mathrm{BV}_\kappa$ such that $q\leq p$ and $M\proj\lambda\in\M_q$.
\end{proposition}
\begin{prooff}
    \item Let $\hh:=(H((2^\kappa)^+),\in,U)$.
    Up to modifying $U$, we may assume that the structure $\hh$ has a lightface definable wellordering and that $\kappa,p,\Hi,\dd$ are lightface definable singeltons over $\hh$.

    \item Let $\lambda<\kappa$ be such that there exists $H\prec\hh$ satisfying that $\kappa\cap H=\lambda$, let $g\in V[h]$ be $V$-generic for $\cc^\mathrm{BV}_\kappa\cap V_\lambda=\cc^\mathrm{BV}_\lambda$ containing $p$, and let $h$ be $V$-generic for $\Col(\omega,<\kappa)$ such that $g\in V[h]$.
    We work in $V[h]$.

    \item Let $F:=\{w_p : p\in g\}$.
    By Lemma \ref{1535}, $F$ is a $\dd$-generic filter for $\Hi$.

    \item Since $(\Hi,\dd)$ is AS-good at $\kappa$ (in $V$), there exist (in $V[h]$) an elementary embedding $\tau:V\to W$ with
    $$\crit(\tau)=\omega_1^V\in\tau(S)$$
    and a $\tau(\dd)$-generic filter $F^*\subseteq\tau(\Hi)$ such that $\tau [F]\subseteq F^*$.

    \item Let $N:=\hull^{\tau(\hh)}(\omega_1^V\cup\{\tau(p)\})\prec\tau(\hh)$.
    It holds that
    \begin{parts}
        \item $N\in W$,

        \item $|N|^W=\omega$,

        \item $N\prec\tau [\hh]$,

        \item $\delta_N=\omega_1^V\in\tau(S)$.
    \end{parts}

    \item\label{465'}\claim In $W$, there exists $q\in\tau(\cc^{\mathrm{BV}}_\kappa)$ such that $q\leq \tau(p)$ and $N\proj\tau(\lambda)\in\M_q$.
    \begin{prooff}
        \item We will verify that one can take
        $$q:=(\tau(w_p),\tau(\M_p)\cup\{N\proj\tau(\lambda)\}).$$
        To that end, it suffices to show that $q\in\tau(\cc^{\mathrm{BV}}_\kappa)$.

        \item Let us assume otherwise.
        Then there exists a winning strategy $\sigma$ for Player I in $\Game^{\mathrm{BV}}_\kappa(q)^W$.
        We will reach a contradiction by defeating this strategy in $V[h]$.
        (We use in this step that the game $\Game^{\mathrm{BV}}_\kappa(q)^W$ is closed for Player II.)

        \item Let Player I play according to $\sigma$ and let us show that Player II can play by maintaining that for all $n\in [-1,\omega)$, there exist $r,X,P$ such that
        \begin{parts}
            \item $r\in g$,

            \item $w_{p_n}\in F^*$ and $w_{p_n}\leq\tau(w_r)$,

            \item $X$ is a finite subset of $\tau [V_\lambda]$,
            
            \item $P=\hull^{\tau(\hh)}(\omega_1^V\cup X)$,       

            \item $\M_{p_n}=\tau(\M_r)\cup\{P\proj\tau(\lambda)\}$.
        \end{parts}
        
        \item The above conditions are satisfied for $n=-1$ and $p_{-1}=q$.
        Let us consider the case $n\geq 0$ and the move $Q_n$ of Player I.
        We let $r,X,P$ be as above, but for the index $n-1$.

        \item\textbf{Case I.} $Q_n=D$ for some $D\in\tau(\dd)$.
        \begin{proof}
            Since $F^*$ is $\tau(\dd)$-generic, there exists $w\in F^*\cap D$.
            Let $w_{p_n}:=w_{p_{n-1}}w\in F^*$ and let $\M_{p_n}:=\M_{p_{n-1}}$.
            It is easily seen that $p_n:=(w_{p_n},\M_{p_n})$ is as required.
        \end{proof}

        \item\textbf{Case IIa.} $Q_n=(P\proj\tau(\lambda),E)$ for some $E\in P\proj\tau(\lambda)$.
        \begin{proof}
            Since $P\subseteq\im(\tau)$, it follows that $P\proj\tau(\lambda)\subseteq\ran(\tau)$.
            Hence, there is $\Bar{E}\in V$ such that $\tau(\Bar{E})=E$.
            By elementarity, $\bar E$ is dense in $\cc^\mathrm{BV}_\lambda$ (in $V$).
            Since $g$ is $V$-generic, there exists $s\in\bar E\cap g$ such that $s\leq r$.
            Observe that $\tau(s)\in\tau(\bar E)=E$ and that
            $$P':=\hull^{\tau(\hh)}(\omega_1^V\cup X\cup\{\tau(s)\})$$
            satisfies that $\delta(P')=\delta(P)=\omega_1^V$.
            Thus, we may set 
            $$p_n:=(\tau(w_s)w_{p_{n-1}},\tau(\M_s)\cup\{P'\proj\tau(\lambda)\}\in\tau(\P^\mathrm{sp}\cap V_\kappa).$$
        \end{proof}



            




        \item\textbf{Case IIb.} $Q_n=(M,E)$ for some $M\in\tau(\M_r)$ and some $E\in M$.
        \begin{proof}
            Observe that the preimages $\bar M:=\tau^{-1}(M)\in\M_r$ and $\bar E:=\tau^{-1}(E)\in\bar M$ are defined.
            By elementarity of $\tau$, $\bar E$ is dense in $\cc^\mathrm{BV}_{\lambda_{\bar M}}$ (in $V$).
            We will verify below that the set
            \begin{equation}\label{246}
                \{s\leq r : \exists t\in\bar E,\, s\leq t,\, \delta(\hull(\bar M,t))=\delta(\bar M)\}
            \end{equation}
            is dense below $r$ (in $\cc^\mathrm{BV}_\lambda$, in $V$).
            Once we know this, the genericity of $g$ yields $s\in g$ and $t\in\bar E$ such that $s\leq r,t$ and $\delta(\hull(\bar M,t))=\delta(\bar M)$.
            For $P':=\hull^{\tau(\hh)}(\omega_1^V\cup X\cup\{\tau(s)\})$, we have that $\delta_{P'}=\omega_1^V=P$.
            Thus, we may set 
            $$p_n:=(\tau(w_s) w_{p_{n-1}},\tau(\M_s)\cup\{P'\proj\tau(\lambda)\})\in\tau(\P^\mathrm{sp}).$$

            \hspace{1em} Let us now verify that the set (\ref{246}) is dense.
            Let $r'\in\cc^\mathrm{BV}_\lambda$ be an arbitrary satisfying that $r'\leq r$ and let us find $s\leq r'$ which belongs to the set (\ref{246}).
            Then there exists $\bar M'\in\M_{r'}$ such that $\delta(\bar M')=\delta(\bar M)$, $\lambda_{\bar M'}=\lambda_{\bar M}$, and $\bar M\subseteq\bar M'$.
            Using the winning strategy for Player II in $\Game^\mathrm{BV}_{\lambda}(r')$, we obtain $s\in\cc^\mathrm{BV}_\lambda$ such that $s\leq r'$ and such that there exists $t\in\bar E$ satisfying $s\leq t$ and $\delta(\hull(\bar M',t))=\delta(\bar M')$.
            We have that
            $$\delta(\bar M)\leq\delta(\hull(\bar M,t))\leq\delta(\hull(\bar M',t))=\delta(\bar M')=\delta(\bar M),$$
            which shows that all of the above values are equal.
            Hence, $s\leq r'$ is as required.
        \end{proof}

        \item Cases I, IIa, and IIb show that Player II can survive $\omega$ steps against the strategy $\sigma$, which is a contradiction.
    \end{prooff}

    \item The conclusion of the proposition now follows from Claim \ref{465'}, by the elementarity of $\tau$.
\end{prooff}

We are now ready to show how Player II can win $\Game^\mathrm{ssp}_\kappa(o)$, provided that $(\Hi,\dd)$ is AS-good at $\kappa$.

\begin{theorem}\label{1696}
    Suppose that
    \begin{assume}
        \item $\kappa>\rank(\Hi)$ is inaccessible,

        \item $(\Hi,\dd)$ is AS-good at $\kappa$.
    \end{assume}
    Then Player II has a winning strategy in $\Game^\mathrm{ssp}_\kappa(o)$.
\end{theorem}
\begin{prooff}
    \item We show that Player II can play by maintaining that for all $n\in [-1,\omega)$, $p_n\in\cc^\mathrm{BV}_\kappa$.

    \item For $n=-1$, Player II can simply play $p_{-1}:=o$.

    \item Let us consider the move $Q_n$ of Player I.

    \item\textbf{Case I.} $Q_n=D$ for some $D\in\dd$.
    \begin{proof}
        Let us consider a partial play
        \begin{center}
            \begin{tabular}{c|ccc}
                I &           & $D$ &      \\
                \hline
                II & $q_{-1}$ &     & $q_0$
            \end{tabular}
        \end{center}
        of $\Game^{\mathrm{BV}}_\kappa(p_{n-1})$ according to some winning strategy for Player II.
        The rules imply that for some $w\in D$, $w_{q_0}\leq w$.
        They also imply that $q_0\leq q_{-1}\leq p_{n-1}$.
        This is enough to conclude that Player II can play $p_n:=q_0$ in the play of $\Game^\mathrm{ssp}_\kappa(o)$.
    \end{proof}

    \item\textbf{Case II.} $Q_n=(U,S)$ for some $U\subseteq H((2^\kappa)^+)$ and for some $S\subseteq\omega_1$ which is stationary.
    \begin{proof}
        By Proposition \ref{134}, there exist a countable 
        $$M\prec (H((2^\kappa)^+,\in,U,\kappa,p_{n-1})$$
        with $\delta_M\in S$, an ordinal $\lambda\in\Bf\cap\kappa$ with $\Hull(M,V_\lambda)\cap\kappa\subseteq\lambda$, and a condition $p_n\in\cc^\mathrm{BV}_\kappa$ such that $p_n\leq p_{n-1}$ and $M\proj\lambda\in\M_{p_n}$.
        We see now that Player II can respond by playing $p_n$.
    \end{proof}

    \item\textbf{Case III.} $Q_n=(M,E)$ for some $M\in\M_{p_{n-1}}$ and some $E\in M$.
    \begin{prooff}
        \item We have that there exists $M^*\prec H((2^\kappa)^+)$ and a lifting
        $$\pi:M\to M^*$$
        such that $p_{n-1}\in \pi(E)$.

        \item Let
        \begin{parts}
            \item $F:=\cc^\mathrm{BV}_{\lambda_M}\cap E$,

            \item $F^\bot:=\{q\in\cc^\mathrm{BV}_{\lambda_M} : \forall r\in F, q\bot_{\cc^\mathrm{BV}_{\lambda_M}} r\}$,

            \item $G:=F\cup F^\bot$.
        \end{parts}
        We have that
        \begin{parts}
            \item $F,F^\bot,G\in M$,

            \item $G$ is dense in $\cc^\mathrm{BV}_{\lambda_M}$,

            \item $p_{n-1}\in\pi(F)$.
        \end{parts}

        \item Let us consider a partial play
        \begin{center}
            \begin{tabular}{c|ccc}
                I &           & $(M,G)$ &      \\
                \hline
                II & $q_{-1}$ &     & $q_0$
            \end{tabular}
        \end{center}
        of $\Game^{\mathrm{BV}}_\kappa(p_{n-1})$ according to some winning strategy for Player II.

        \item The rules imply that there exists $r\in G$ such that $\delta(\Hull(M,r))=\delta(M)$ and that $q_0\leq r,q_{-1}$.

        \item Note that $r\in V_{\lambda_M}$.
        This means that $r=\pi(r)\in \pi(F)\cup\pi(F^\bot)$.

        \item Since $q_0\leq r$ and $q_0\leq q_{-1}\leq p_{n-1}$, we see that $r\para p_{n-1}$ in the sense of $\cc^\mathrm{BV}_{\lambda_M}$.
        However, this also means that $r\para p_{n-1}$ in the sense of $\cc^\mathrm{BV}_{\kappa}$.

        \item Since $p_{n-1}\in\pi(F)$, we get that $r\not\in\pi(F^\bot)$.
        This means that $r\in\pi(F)$ and consequently, $r\in F$.

        \item Since $F\subseteq E$,  we get that
        \begin{parts}
            \item $r\in E$ and $\delta(\Hull(M,r))=\delta(M)$,

            \item $q_0\leq p_{n-1},r$,

            \item $q_0\in\cc^\mathrm{BV}_{\kappa}$.
        \end{parts}
        In other words, Player II can respond in the play of $\Game^\mathrm{ssp}_\kappa(o)$ by playing $p_n:=q_0$.
    \end{prooff}

    \item The verification of three cases shows that Player II can always play $\omega$ many rounds of $\Game^\mathrm{ssp}_\kappa(o)$ without breaking the rules.
    This concludes the argument since all such plays are won by Player II.
\end{prooff}

\begin{corollary}
    Suppose that
    \begin{assume}
        \item $\kappa>\rank(\Hi)$ is inaccessible,

        \item $(\Hi,\dd)$ is AS-good at $\kappa$.
    \end{assume}
    Then $\cc^\mathrm{ssp}_\kappa$ is a stationary set preserving poset that adds a $\dd$-generic filter for $\Hi$.
    \qed
\end{corollary}

Our arguments actually show that $\cc^\mathrm{BV}_\kappa\subseteq\cc^\mathrm{ssp}_\kappa$ and that this set is closed for some winning strategy for Player II in $\Game^\mathrm{ssp}_\kappa(-)$.
The point is that the winning strategy for Player II in $\Game^\mathrm{ssp}_\kappa(-)$ is simply ``play so as to remain inside $\cc^\mathrm{BV}_\kappa$''.
This ensures that $\cc^\mathrm{BV}_\kappa$ is yet another ssp poset adding a $\dd$-generic filter for $\Hi$.
We will have to come back to this poset in Section \ref{An Example for the AS Goodness}, so we observe the following corollary as well.

\begin{corollary}
    Suppose that
    \begin{assume}
        \item $\kappa>\rank(\Hi)$ is inaccessible,

        \item $(\Hi,\dd)$ is AS-good at $\kappa$.
    \end{assume}
    Then $\cc^\mathrm{BV}_\kappa$ is an initial segment of $\cc^\mathrm{ssp}_\kappa$, it is stationary set preserving, and it adds a $\dd$-generic filter for $\Hi$.
\end{corollary}
\begin{proof}
    The proof of Theorem \ref{1696} actually shows that for all $p\in\cc^\mathrm{BV}_\kappa$, Player II has a winning strategy $\sigma$ in $\Game^\mathrm{ssp}_\kappa(p)$ such that $\ran(\sigma)\subseteq\cc^\mathrm{BV}_\kappa$.
    This show that $\cc^\mathrm{BV}_\kappa\subseteq\cc^\mathrm{ssp}_\kappa$ and that $\cc^\mathrm{BV}_\kappa$ is a stationary set preserving poset adding a $\dd$-generic filter for $\Hi$ (Lemma \ref{1806}).
    We observed in Lemma \ref{1535} that $\cc^\mathrm{BV}_\kappa$ is an initial segment of $\P^\mathrm{sp}\cap V_\kappa$, whence the final conclusion.
\end{proof}

\begin{example}\label{2398}
    Suppose that $\ns$ is saturated and that $\kappa$ is an inaccessible cardinal.
    We want to show that there exists a stationary set preserving poset adding a cofinal function $f:\omega\to\omega_2^V$.
    This comes down to showing that the Namba problem $(\Hi,\dd)$ of Example \ref{Namba problem} is spp-consistent.
    A sufficient condition for the ssp consistency is the AS goodness at $\kappa$ and this problem satisfies that.
    More precisely, let $S\subseteq\omega_1$ be stationary, let $g$ be $V$-generic for $\Col(\omega,<\kappa)$, and let $F$ be a $\dd$-generic filter for $\Hi$.
    We want to find
    $$\tau:V\to W$$
    satisfying that $\crit(\tau)=\omega_1^V\in\tau(S)$ and a $\tau(\dd)$-generic filter $F^*$ for $\tau(\Hi)$ satisfying that $\tau [F]\subseteq F^*$.
    
    Let $G$ be a $V$-generic filter for $(\ns^+)^V$ containing $S$ and let
    $$\tau:V\to W$$
    be the corresponding generic embedding.
    Observe that $f:=\bigcup F$ is a cofinal function $\omega\to\omega_2^V$.
    Since $\tau$ is continuous at $\omega_2^V$, it follows that $\tau\circ f$ is a cofinal function $\omega\to\omega_2^W$.
    Hence,
    $$F^*:=\{w\in\tau(\Hi) : w\subseteq\tau\circ f\}$$
    is a $\tau(\dd)$-generic filter for $\tau(\Hi)$.
    This shows that $\tau$ and $F^*$ are as required.
    
    We have shown that the Namba problem is AS-good at $\kappa$.
    This means that there exists a stationary set preserving poset adding a solution to it.
    This solution is exactly a cofinal function $\omega\to\omega_2^V$.
    Furthermore, we know that $\cc^\mathrm{ssp}_\kappa$ and $\cc^\mathrm{BV}_\kappa$ are stationary set preserving posets adding such a function.
    \qed
\end{example}


\section{An Example for the AS Goodness}
\label{An Example for the AS Goodness}

In this section, we show how our methods can be used to give another proof of the main lemma of \cite{aspero2021martin}.
Before we precisely state what we are doing, we need a following auxiliary notion.

\begin{definition}
    Suppose that $\theta$ is a cardinal and $\Gamma\subseteq\ps(\R)$.
    Then $\Gamma$ is \intro{$\theta$-well behaved} iff
    \begin{parts}
        \item for all $A\in\Gamma$, $A$ is $\theta$-universally Baire,

        \item for all $A\in\Gamma$, for all formulas $\phi(x,\dot A)$ projective in $\dot A$, we have 
        $$B:=\{x\in\R : \phi(x,A)\}\in\Gamma$$
        and $V^{\Col(\omega,\theta)}\models B_\G=\{x\in\R_{\G} : \phi(x,A_\G)\}$.
        \qed
    \end{parts}
\end{definition}

\begin{declaration}
    Suppose that
    \begin{assume}
        \item $\ns$ is saturated,

        \item $2^{\omega_1}=\omega_2$,

        \item $\mathsf{MA}_{\omega_1}$ holds,

        \item $D$ is dense in $\pmax$,

        \item for all $\theta$, there exists $\theta$-well behaved $\Gamma$ such that $D$ is coded by a set in $\Gamma$,

        \item $A\subseteq\omega_1$ satisfies that $\omega_1^{L[A]}=\omega_1$,
        
        \item $\kappa$ is inaccessible.
        \qed
    \end{assume}
\end{declaration}

Our goal is to construct an ssp poset $\Q$ such that in $V^\Q$, there exist generic iterations $\mathcal{I}$ and $\mathcal{J}$ of $\pmax$ conditions of length $\omega_1+1$ such that
\begin{parts}
    \item $M_0^\mathcal{J}\in D_{\G}$,
    
    \item $\mathcal{I}\rest\omega_1^{M_0^\mathcal{J}}$ witnesses that $M_0^\mathcal{J}<M_0^\mathcal{I}$, 

    \item $\pi^{\mathcal{J}}(\mathcal{I}\rest\omega_1^{M_0^\mathcal{J}})=\mathcal{I}$,

    \item $M_{\omega_1}^\mathcal{I}=(H_{\omega_2},\in,\ns,A)^V$,

    \item\label{635} $\ns\cap\ps(\omega_1)^{M^\jj_{\omega_1}}=I^{M^\jj_{\omega_1}}$.
\end{parts}
We will start by forgetting for a moment about the requirement \ref{635} and recasting the rest of the problem as adding a model to an infinitary propositional formula by an ssp forcing.
It will turn out that we can get \ref{635} for free.

        





\begin{notation}
    Suppose $\M=(\omega,\dots)$ is a structure in a recursively enumerable language.
    Then \intro{$\vartheta(\M)$} denotes the truth function of the model $(\M,\, m : m<\omega)$, i.e. $\vartheta(\M):\omega\to 2$ is defined in such a way that for all $n<\omega$, for $\psi$ denoting the formula in the language of
    $$(\M,\, m : m<\omega)$$
    with $\ulcorner\psi\urcorner=n$, it holds that 
    $$\vartheta(\M)(n)=1\iff (\M,\, m : m<\omega)\models\psi.$$
    \qed
\end{notation}

\begin{notation}
\hfill
    \begin{parts}
        \item If $E\subseteq X\times Y$ and $x\in X$, we denote $E_x:=\{y\in Y : (x,y)\in E\}$.

        \item If $a:X\to Y$ and $x\in X$, then $a_x:=a(x)$.\qed
    \end{parts}
\end{notation}

It is not clear that we can express the fact that something is a $\pmax$ condition using propositional logic.
Our workaround is to fix a tree that projects to codes of $\pmax$ conditions and just make sure that a code of the object that we are forcing belongs to the projection of that tree.
Similarly, we will not express directly that one of the conditions that we are forcing belongs to the reinterpretation of the set $D$, but rather that its code belongs to the projection of the corresponding tree.

\begin{declaration}
    Let $T,U$ be $\omega_2$-absolutely complemented trees on $\omega\times\omega_2$ such that
    \begin{parts}
        \item $p[T]$ is the set of all $\vartheta(\omega,E,I,a)$ for which there exists $p\in\pmax$ satisfying $(\omega,E,I,a)\simeq p$,

        \item $p[U]$ is the set of all $\vartheta(\omega,E,I,a)$ for which there exists $p\in D$ satisfying $(\omega,E,I,a)\simeq p$.\qed
    \end{parts}
\end{declaration}

Having made the necessary compromises, we are ready to restate the goal in a way that can easily be coded by an infinitary propositional formula.

\begin{goal}
\label{goal reformulated}
    Construct an ssp poset $\Q$ such that in $V^\Q$, there exist $E$, $F$, $I$, $J$, $a$, $b$, $u$, $v$, $\pi$, $\sigma$, $\zh$, $r$ such that
    \begin{parts}
        \item\label{220} $E,F\subseteq\omega_1\times(\omega\times\omega)$,

        \item $I,J\subseteq\omega_1\times\omega$,

        \item $a,b:\omega_1\to\omega$,

        \item $u,v:\omega\to\omega_2$,

        \item\label{465} for all $n<\omega$, $(\vartheta(\omega,E_0,I_0,a_0)\rest n,u\rest n)\in T$,

        \item\label{467} for all $n<\omega$, $(\vartheta(\omega,F_0,J_0,b_0)\rest n, v\rest n)\in U$,

        \item $\pi,\sigma\subseteq\{(\alpha,\beta) : \alpha\leq\beta<\omega_1\}\times(\omega\times\omega)$,

        \item for all $\alpha\leq\beta<\omega_1$, $\pi_{\alpha,\beta}:(\omega,E_\alpha,I_\alpha,a_\alpha)\xrightarrow[\Sigma_\omega]{}(\omega,E_\beta,I_\beta,a_\beta)$,

        \item for all $\alpha\leq\beta\leq\gamma<\omega_1$, $\pi_{\beta,\gamma}\circ\pi_{\alpha,\beta}=\pi_{\alpha,\gamma}$,

        \item for all $\alpha<\omega_1$, there exists $k<\omega$ such that
        \begin{parts}
            \item $(\omega,E_{\alpha+1})\models\dot k<\dot\omega_1$,

            \item for all $\ch<\omega$, there exists $\dzh<\omega$ such that $(\omega,E_\alpha)\models\dot \dzh:\dot\omega_1\to\dot V$ and such that for $\sh:=\pi_{\alpha,\alpha+1}(\dzh)$, $(\omega,E_{\alpha+1})\models \dot\sh(k)=\dot\ch$,

            \item for $g_\alpha:=\{w<\omega : (\omega,E_\alpha,I_\alpha)\models \dot w\in I_\alpha^+,\ k\, E_{\alpha+1}\,\pi_{\alpha,\alpha+1}(w)\}$, we have that $(\omega,E_\alpha,I_\alpha,g_\alpha)\models``g_\alpha$ is generic for $I_\alpha^+"$,
        \end{parts}

        \item for all limit $\gamma<\omega_1$, for all $n<\omega$, there exist $\alpha<\gamma$ and $m<\omega$ such that $\pi_{\alpha,\gamma}(m)=n$,

        \item for all $\alpha\leq\beta<\omega_1$, $\sigma_{\alpha,\beta}:(\omega,F_\alpha,J_\alpha,b_\alpha)\xrightarrow[\Sigma_\omega]{}(\omega,F_\beta,J_\beta,b_\beta)$,

        \item for all $\alpha\leq\beta\leq\gamma<\omega_1$, $\sigma_{\beta,\gamma}\circ\sigma_{\alpha,\beta}=\sigma_{\alpha,\gamma}$,

        \item for all $\alpha<\omega_1$, there exists $k<\omega$ such that
        \begin{parts}
            \item $(\omega,F_{\alpha+1})\models\dot k<\dot\omega_1$,

            \item for all $\ch<\omega$, there exists $\dzh<\omega$ such that $(\omega,F_\alpha)\models\dot\dzh:\dot\omega_1\to\dot V$ and such that for $\sh:=\sigma_{\alpha,\alpha+1}(\dzh)$, $(\omega,F_{\alpha+1})\models \dot\sh(k)=\dot\ch$,

            \item for $h_\alpha:=\{w<\omega : (\omega,F_\alpha,J_\alpha)\models \dot w\in J_\alpha^+,\ k\, F_{\alpha+1}\,\sigma_{\alpha,\alpha+1}(w)\}$, we have that $(\omega,F_\alpha,J_\alpha,h_\alpha)\models``h_\alpha$ is generic for $J_\alpha^+"$,
        \end{parts}

        \item for all limit $\gamma<\omega_1$, for all $n<\omega$, there exist $\alpha<\gamma$ and $m<\omega$ such that $\sigma_{\alpha,\gamma}(m)=n$,

        \item $\zh\subseteq\omega_1\times(\omega\times H_{\omega_2}^V)$,

        \item for all $\alpha<\omega_1$, $\zh_\alpha:(\omega,E_\alpha,I_\alpha,a_\alpha)\xrightarrow[\Sigma_\omega]{}(H_{\omega_2},\in,\ns,A)^V$,

        \item for all $\alpha\leq\beta<\omega_1$, $\zh_\alpha=\zh_\beta\circ\pi_{\alpha,\beta}$,

        \item\label{270} for all $x\in H_{\omega_2}^V$, there exist $\alpha<\omega_1$ and $\ch<\omega$ such that $\zh_\alpha(\ch)=x$,

        \item $r:\omega_1\to(\omega\to\omega_1)$,

        \item\label{296} for all $\alpha<\omega_1$, $r_\alpha$ is the ranking function of $(\omega,F_\alpha)$.\qed
    \end{parts}
\end{goal}

Note that if $|\Q|>\omega_2$, then every element of $p[T]$ in $V^\Q$ will be of the form $\vartheta(\M)$ for some $\M$ isomorphic to a $\pmax$-condition, but the converse might not hold, and similarly for $U$.
However, this suffices for our purposes.
Here are some examples how to convert the requirements of Goal \ref{goal reformulated} into a propositional formula.

\begin{parts}
    \item Let $\dot E_{\alpha,m,n},\dot F_{\alpha,m,n}, \dot I_{\alpha,m}, \dot J_{\alpha,m},\dot a_{\alpha,m},\dot b_{\alpha,m},\dot u_{m,\xi},\dot v_{m,\xi},\dot\pi_{\alpha,\beta,m,n}$, $\dot\sigma_{\alpha,\beta,m,n}$,\\
    $\dot\zh_{\alpha,m,x}$, $\dot r_{\alpha,m,\gamma}$ be propositional letters, for all $\alpha\leq\beta<\omega_1$, $\gamma<\omega_1$, $m,n<\omega$, $\xi<\omega_2$, $x\in H_{\omega_2}$.

    \item If $\psi$ is the language of some $(\omega,E_\alpha,I_\alpha,a_\alpha)$, we can express
    $$(\omega,E_\alpha,I_\alpha,a_\alpha)\models\psi$$
    as a propositional formula by induction on complexity of $\psi$.
    Similarly for $(\omega,F_\alpha,J_\alpha,b_\alpha)$.

    \item We can now express requirement \ref{465} as
    $$\bigwedge_{n<\omega}\bigvee_{\substack{s\in T,\\ |s|=n}}\bigwedge_{i<n}\theta_{s,i}$$
    where
    \begin{parts}
        \item for $s(i)=(1,\xi)$ and for $\psi_i$ being the formula with code $i$, $\theta_{s,i}$ is defined as the propositional formula asserting
        $$``(\omega,E_\alpha,I_\alpha,a_\alpha)\models\psi_i"\wedge \dot u_{i,\xi},$$

        \item for $s(i)=(0,\xi)$ and for $\psi_i$ being the formula with code $i$, $\theta_{s,i}$ is defined as the propositional formula asserting
        $$``(\omega,E_\alpha,I_\alpha,a_\alpha)\models\neg\psi_i"\wedge \dot u_{i,\xi}.$$
    \end{parts}
    Similarly for the requirement \ref{467}.
\end{parts}

One can proceed like this to code all of the requirements of Goal \ref{goal reformulated} and obtain a formula which will have $\omega_2$ propositional letters and whose every conjunction and disjunction will have the length at most $\omega_2$.

\begin{definition}
\label{548}
    We denote by \intro{$\phi_\mathrm{AS}$} the negation normal form of the formula sketched just above.\qed
\end{definition}

The poset $\Q$ required in Goal \ref{goal reformulated} will be the poset $\cc^\mathrm{BV}_\kappa$ of Definition \ref{1488}, defined with respect to the problem $(\Hi(\phi_\mathrm{AS}),\dd(\phi_\mathrm{AS}))$.
The several things need to be verified, starting from the fact that the said problem is well-defined.

\begin{lemma}
    $\phi_\mathrm{AS}$ is consistent.
\end{lemma}
\begin{proof}
    This is shown as \cite[Lemma 3.6]{aspero2021martin}.
\end{proof}

The following lemma simply states that the formula $\phi_\mathrm{AS}$ codes that which it is designed to code.

\begin{lemma}
    Suppose that $\mu\models\phi_\mathrm{AS}$ belongs to an outer model.
    Then in the outer model, there exist unique generic iterations $\mathcal{I}$ and $\mathcal{J}$ of $\pmax$ conditions of length $\omega_1^V+1$ such that
    \begin{parts}
        \item $\mathcal{I}\rest\omega_1^{M_0^\mathcal{J}}$ witnesses that $M_0^\mathcal{J}<M_0^\mathcal{I}$, 

        \item $\pi^{\mathcal{J}}(\mathcal{I}\rest\omega_1^{M_0^\mathcal{J}})=\mathcal{I}$,

        \item $M_{\omega_1}^\mathcal{I}=(H_{\omega_2},\in,\ns,A)^V$,

        \item $M_0^\mathcal{J}\in D_g$,

        \item for all $\alpha<\omega_1^V$, $(\omega,E_\alpha,I_\alpha,a_\alpha)\simeq M_\alpha^\ii$ and $(\omega,F_\alpha,J_\alpha,b_\alpha)\simeq M_\alpha^\jj$,

        \item for all $\alpha\leq\beta<\omega_1^V$, the diagrams
        \begin{center}
            \begin{tikzcd}
                (\omega,E_\alpha,I_\alpha,a_\alpha)\arrow[r, "\pi_{\alpha,\beta}"]\arrow[d, "\simeq"] & (\omega,E_\beta,I_\beta,a_\beta)\arrow[d, "\simeq"]\\
                M_\alpha^\ii\arrow[r, "\pi^\ii_{\alpha,\beta}"] & M_\beta^\ii
            \end{tikzcd}
        \end{center}
        \begin{center}
            \begin{tikzcd}
                (\omega,E_\alpha,I_\alpha,a_\alpha)\arrow[r, "\zh_{\alpha}"]\arrow[d, "\simeq"] & (H_{\omega_2},\in,\ns,A)^V\arrow[d, "="]\\
                M_\alpha^\ii\arrow[r, "\pi^\ii_{\alpha,\omega_1^V}"] & (H_{\omega_2},\in,\ns,A)^V
            \end{tikzcd}
        \end{center}
        \begin{center}
            \begin{tikzcd}
                (\omega,F_\alpha,J_\alpha,b_\alpha)\arrow[r, "\sigma_{\alpha,\beta}"]\arrow[d, "\simeq"] & (\omega,F_\beta,J_\beta,b_\beta)\arrow[d, "\simeq"]\\
                M_\alpha^\jj\arrow[r, "\pi^\jj_{\alpha,\beta}"] & M_\beta^\jj
            \end{tikzcd}
        \end{center}
        commute.\qed
    \end{parts}
\end{lemma}

\begin{notation}
    In the case of the above proposition, we denote by \intro{$\ii^\mu$} and \intro{$\jj^\mu$} the generic iterations $\ii$ and $\jj$, respectively.
    We also denote 
    $$E^\mu:=\{(\alpha,m,n)\in \omega_1^V\times(\omega\times\omega) : \mu(\dot E_{\alpha,m,n})=1\},$$
    $$F^\mu:=\{(\alpha,m,n)\in \omega_1^V\times(\omega\times\omega) : \mu(\dot F_{\alpha,m,n})=1\},$$
    $$I^\mu:=\{(\alpha,m)\in\omega_1^V\times\omega : \mu(\dot I_{\alpha,m})=1\},$$
    and so on.\qed
\end{notation}

We now want to show that $\phi_\mathrm{AS}$ is AS-good at $\kappa$.
By this, we simply mean that the problem $(\Hi(\phi_\mathrm{AS}),\dd(\phi_\mathrm{AS}))$ is AS-good at $\kappa$.
After unraveling the definitions, it turns out that we need to verify that for all stationary $S\subseteq\omega_1$, it holds in $V^{\Col(\omega,<\kappa)}$ that for all $\mu\models\phi$, there exist elementary $\tau:V\to W$ and $\hat\mu\models\tau(\phi)$ such that $\crit(\tau)=\omega_1^V\in\tau(S)$ and for all subformulas $\psi$ of $\phi$, $\hat\mu(\tau(\psi))=\mu(\psi)$.

\begin{proposition}
\label{827}
    $\phi_\mathrm{AS}$ is AS-good at $\kappa$.
\end{proposition}
\begin{prooff}
    \item Let $h$ be $V$-generic for $\Col(\omega,<\kappa)$ and let us work in $V[h]$.
    We fix an arbitrary valuation $\mu\models\phi$ and an arbitrary stationary $S\subseteq\omega_1^V$ and we are going to produce elementary $\tau:V\to W$ and $\hat\mu\models\tau(\phi)$ such that $\crit(\tau)=\omega_1^V\in\tau(S)$ and for all subformulas $\psi$ of $\phi$, $\hat\mu(\tau(\psi))=\mu(\psi)$.
    
    \item Let $\hat\jj$ be a generic iteration of length $\omega_1+1$ extending $\jj^\mu$ and satisfying $S\in g_{\omega_1^V}^{\hat\jj}$ and let $\hat\ii:=\pi^{\hat\jj}_{\omega_1^V,\omega_1}(\ii)$.
    Since 
    $$M^{\ii}_{\omega_1^V}=(H_{\omega_2},\in,\ns,A)^V,$$
    there exists a unique generic iteration $\kk$ of $V$ of length $\omega_1+1$ such that for all $\alpha<\omega_1$, $g^\kk_\alpha=g^{\hat\ii}_{\omega_1^V+\alpha}$.

    \item We let $\tau:=\pi^\kk:V\to M^\kk_{\omega_1}=:W$.

    \item For all $\alpha<\omega_1$, let $e_\alpha:\omega\to M^{\hat\ii}_\alpha$ be a bijection such that if $\alpha<\omega_1^V$, then $e_\alpha$ is the transitive collapse of $(\omega,E_\alpha)$.
    We define for all $\alpha<\omega_1$,
    $$\hat E_\alpha:=\{(m,n) : e_\alpha(m)\in e_\alpha(n)\}.$$
    In the case $\alpha<\omega_1^V$, we get the old $\hat E_\alpha=E_\alpha$.

    \item In the similar fashion, we extended other sets from Goal \ref{goal reformulated}, producing the ``hat versions'' of them.
    Then the valuation $\hat\mu$ is defined according to these ``hat versions''.

    \item\claim For $\psi$ a subformula of $\phi$, $\hat\mu(\tau(\psi))=\mu(\psi)$.
    \begin{proof}
        This requires a tedious verification.
        Suppose for example that $\psi$ comes from the part \ref{270}. of Goal \ref{goal reformulated}.
        Then $\psi$ might of the formula
        $$\bigwedge_{x\in H_{\omega_2}^V}\bigvee_{\alpha<\omega_1^V}\bigvee_{\ch<\omega}\dot{\zh}_{\alpha,\ch,x}.$$
        In this case, $\mu(\psi)=1$, $\tau(\psi)$ is the formula
        $$\bigwedge_{x\in H_{\omega_2}^W}\bigvee_{\alpha<\omega_1}\bigvee_{\ch<\omega}\dot{\zh}_{\alpha,\ch,x},$$
        and we need to verify that $\tau(\psi)$ is evaluated as true by $\hat\mu$.
        However, this is trivial, since 
        $$(H_{\omega_2}^W,\zh_\alpha : \alpha<\omega_1)$$
        is the direct limit 
        $$(M_\alpha^{\hat\ii},\pi^{\hat\ii}_{\alpha,\beta} : \alpha\leq\beta<\omega_1).$$
        We leave it to the reader to convince herself of the rest.
    \end{proof}

    \item Thus, $\tau$ and $\hat\mu$ are as required.
\end{prooff}

We now know that both $\cc^\mathrm{ssp}_\kappa$ and $\cc^\mathrm{BV}_\kappa$ are stationary set preserving and add a model for $\phi_\mathrm{AS}$.
Let $g$ be a $V$-generic for one of them and let $\mu_g:=\mu$, $\ii:=\ii^\mu$, and $\jj:=\jj^\mu$.
The $\pmax$-iterations $\ii$ and $\jj$ are like what we wanted, except the thing that we forgot about.
Namely, recall that we did not code into $\phi_\mathrm{AS}$ the requirement that
$$\ns\cap\ps(\omega_1)^{M^{\jj}_{\omega_1}}=I^{M^{\jj}_{\omega_1}}.$$
It turns out that if we force with the smaller of the posets, i.e. $\cc^\mathrm{BV}_\kappa$, then we get this property for ``free''.

\begin{proposition}
    It holds in $V^{\cc^\mathrm{BV}_\kappa}$ that
    $$\ns\cap\ps(\omega_1)^{M^{\jj^{\mu_{\dot g}}}_{\omega_1}}=I^{M^{\jj^{\mu_{\dot g}}}_{\omega_1}}.$$
\end{proposition}
\begin{prooff}
    \item Let us assume otherwise.
    Since ($\supseteq$) clearly holds, the other inclusion must fail.
    More precisely, there exist $p\in\cc^\mathrm{BV}_\kappa$, $\dot C\in V^{\cc^\mathrm{BV}_\kappa}$, $\alpha<\omega_1$, and $m<\omega$ such that $p$ forces that
    \begin{parts}[ref=\arabic{pfenumi}$^\circ$\alph{partsi}]
        \item $(\omega,F^{\mu_{\dot g}}_\alpha)\models``\dot m$ is a subset of $\omega_1$ which is not in $J^{\mu_{\dot g}}_\alpha$'',

        \item $\dot C$ is a club in $\omega_1^V$,

        \item\label{886} the image of $m$ in the final iterate of $\jj^{\mu_{\dot g}}$ has the empty intersection with $\dot C\subseteq V_\kappa$.
    \end{parts}

    \item Let 
    \begin{enumerate}

        \item $\theta:=(2^\kappa)^+$,

        \item $R$ be a wellordering of $H_\theta$,

        \item $\hh:=(H_\theta,\in,R,\kappa,\dot C, T,U,A)$,

        \item $\lambda<\kappa$ be such that there exists $H\prec\hh$ satisfying that $p\in H$ and $\lambda=\kappa\cap H$,
        



        \item $g$ be $V$-generic for $\cc^\mathrm{BV}_\lambda$ containing $p$,

        \item $h$ be $V$-generic for $\Col(\omega,<\kappa)$ such that $g\in V[h]$,

        \item $\mu:=\mu_g\models\phi$.
    \end{enumerate}
    We work in $V[h]$.

    \item Let $S_0$ be the image of $m$ under the transitive collapse of $(\omega,F_\alpha^\mu)$ and let $S:=\pi^{\jj^\mu}_{\alpha,\omega_1^V}(S_0)$.

    \item Let $\hat\ii$, $\hat\jj$, $\kk$, $\tau$, $W$, and $\hat\mu$ be obtained as in the proof of Proposition \ref{827}.
    (The reason why we cannot simply quote this proposition is because now $S$ is not necessarily in $V$.)
    We have that
    \begin{parts}
        \item $\tau:V\to W$ elementarily,
        
        \item $\crit(\tau)=\omega_1^V$,

        \item $\omega_1^V\in\pi^{\hat\jj}_{\omega_1^V,\omega_1}(S)=\pi^{\hat\jj}_{\alpha,\omega_1}(S_0)$,

        \item for all subformulas $\psi$ of $\phi$, $\hat\mu(\tau(\psi))=\mu(\psi)$.
    \end{parts}

    \item Let 
    \begin{assume}
        \item $\djmath_\delta$ be the formula coding the assertion that ``$\omega_1^{M^{\hat\jj}_{\delta}}=\delta$ and this ordinal belongs to the image of the set coded by $m$ under the $\alpha$-to-$(\delta+1)$ iteration embedding'',

        \item $N:=\hull^{\tau(\hh)}(\omega_1^V\cup\{\tau(p)\})$,

        \item $q:=(\tau(w_p)\cup\{\djmath_{\omega_1^V}\},\tau(\M_p)\cup\{N\proj\tau(\lambda)\})\in\tau(\P^\mathrm{sp})$.
    \end{assume}
    
    \item\claim $q\in\tau(\cc^\mathrm{BV}_\kappa)$
    \begin{proof}
        This is shown like Claim \ref{465'} of the proof of Proposition \ref{134}.
    \end{proof}

    \item By elementarity, it holds in $V$ that there exist $r\in\P_\lambda$, $\delta<\omega_1$, and $P\prec\hh$ such that
    \begin{parts}
        \item $r=(w_p\cup\{\djmath_\delta\},\M_p\cup\{P\proj\lambda\})\in\P_\lambda$,

        \item $p\in P$,

        \item $\delta_P=\delta>\alpha$.
    \end{parts}

    \item\claim $r\Vdash^V_{\cc^\mathrm{BV}_\kappa}$``$\delta$ belongs to both the image of the set coded by $m$ along the $\jj$-iteration embedding and to the club $\dot C$''.
    \begin{prooff}
        \item Since $\djmath_\delta\in w_r$, it is immediate that $r\Vdash$``$\delta$ belongs to the image of the set coded by $m$ along the $\jj$-iteration embedding''.
        It remains to verify that $r\Vdash\delta\in\dot C$.

        \item The proof of Claim \ref{1173} of Proposition \ref{1146} shows that $r$ is semigeneric for $(P,\cc^\mathrm{BV}_\kappa)$ (in $V$).
        The point here is that $P\proj\lambda\in\M_r$, so the rules of the game $\Game^\mathrm{BV}_\kappa(r)$ suffice to show the semigenericity.

        \item Let $g'$ be an arbitrary $V$-generic for $\cc^\mathrm{BV}_\kappa$ containing $r$.
        We have that $\delta(P[g'])=\delta(P)=\delta$.

        \item Since $\dot C^{g'}\in P[g']$ is a club in $\omega_1$, we have that $\delta\in C$, as required.
    \end{prooff}

    \item The previous claim is in contradiction with \ref{886}.
\end{prooff}

Finally, we arrive at the conclusion.

\begin{corollary}
    The poset $\cc^\mathrm{BV}_\kappa$ is stationary set preserving and there exist in $V^{\cc^\mathrm{BV}_\kappa}$ generic iterations $\mathcal{I}$ and $\mathcal{J}$ of $\pmax$ conditions of length $\omega_1+1$ such that
    \begin{parts}
        \item $M_0^\mathcal{J}\in D_{\G}$,
        
        \item $\mathcal{I}\rest\omega_1^{M_0^\mathcal{J}}$ witnesses that $M_0^\mathcal{J}<M_0^\mathcal{I}$, 

        \item $\pi^{\mathcal{J}}(\mathcal{I}\rest\omega_1^{M_0^\mathcal{J}})=\mathcal{I}$,

        \item $M_{\omega_1}^\mathcal{I}=(H_{\omega_2},\in,\ns,A)^V$,

        \item $\ns\cap\ps(\omega_1)^{M^\jj_{\omega_1}}=I^{M^\jj_{\omega_1}}$.\qed 
    \end{parts}
\end{corollary}

\printindex
\newpage
\bibliographystyle{alpha}
\bibliography{lit.bib}

\begin{thebibliography}{ARS85}

\bibitem[ACS07]{abraham2007some}
Uri Abraham, James Cummings, and Clifford Smyth.
\newblock Some results in polychromatic {Ramsey} theory.
\newblock {\em J. Symb. Log.}, 72(3):865--896, 2007.

\bibitem[ARS85]{abraham1985on}
Uri Abraham, Matatyahu Rubin, and Saharon Shelah.
\newblock On the consistency of some partition theorems for continuous colorings, and the structure of {{\(\aleph _ 1\)}}-dense real order types.
\newblock {\em Ann. Pure Appl. Logic}, 29:123--206, 1985.

\bibitem[AS21]{aspero2021martin}
David Asper{\'o} and Ralf Schindler.
\newblock {Martin's Maximum{{\(^{++}\)}}} implies {Woodin's} axiom {{\((*)\)}}.
\newblock {\em Ann. Math. (2)}, 193(3):793--835, 2021.

\bibitem[DB24]{bondt2024stationary}
Ben De~Bondt.
\newblock {\em Stationary set preserving forcing, side conditions and games}.
\newblock PhD thesis, Université Paris Cité, 2024.

\bibitem[Far21]{farah2021extender}
Ilijas Farah.
\newblock The extender algebra and {{\(\Sigma^2_1\)}}-absoluteness.
\newblock In {\em Large cardinals, determinacy and other topics. The Cabal Seminar, Vol. IV. Reprints of papers and new material based on the Los Angeles Caltech-UCLA Logic Cabal Seminar 1976--1985}, pages 141--176. Cambridge: Cambridge University Press; Ithaca, NY: Association of Symbolic Logic (ASL), 2021.

\bibitem[FJZ03]{feng2003structure}
Q.~Feng, T.~Jech, and J.~Zapletal.
\newblock On the structure of stationary sets.
\newblock {\em arXiv:math/0311514}, 2003.

\bibitem[KV23]{kasum2023iterating}
Obrad Kasum and Boban Veličković.
\newblock Iterating semi-proper forcing using virtual models.
\newblock {\em arXiv:2303.12565}, 2023.

\bibitem[Lie23]{lietz2023forcing}
Andreas Lietz.
\newblock {\em Forcing {``$\ns$ is $\omega_1$-dense''} from Large Cardinals}.
\newblock PhD thesis, Westf\" alischen Wilhelms-Universit\" at M\" unster, 2023.

\bibitem[Moo06]{moore2006five}
Justin~Tatch Moore.
\newblock A five element basis for the uncountable linear orders.
\newblock {\em Ann. Math. (2)}, 163(2):669--688, 2006.

\bibitem[V{\"a}{\"a}11]{vaananen2011models}
Jouko V{\"a}{\"a}n{\"a}nen.
\newblock {\em Models and games.}, volume 132 of {\em Camb. Stud. Adv. Math.}
\newblock Cambridge: Cambridge University Press, 2011.

\end{thebibliography}
\end{document}